\documentclass[preprint,12pt]{elsarticle}

\usepackage{amssymb}
\usepackage{amsthm}
\usepackage{amsmath}
\usepackage{url}
\newtheorem{thm}{Theorem}
\newtheorem{dfn}[thm]{Definition}
\newtheorem{clr}[thm]{Corollary}
\newproof{pf}{Proof}

\begin{document}

\begin{frontmatter}



\title{Asymptotic approximation of the eigenvalues and the eigenfunctions for the Orr-Sommerfeld equation on infinite intervals}


\author{Victor Nijimbere}

\address{Carleton University, Ottawa, Ontario, Canada}

\begin{abstract}

Asymptotic eigenvalues and eigenfunctions for the Orr-Sommerfeld equation in two and three dimensional incompressible flows on an infinite domain and on a semi-infinite domain are obtained. Two configurations are considered, one in which a short-wave limit approximation is used, and another in which a long-wave limit approximation is used. In the short-wave limit, WKB methods are utilized to estimate the eigenvalues, and the eigenfunctions are approximated in terms Green's functions. 
The procedure consists of transforming the Orr-Sommerfeld equation into a system of two second order ordinary differential equations for which eigenvalues and eigenfunctions can be approximated. The approximated eigenvalues can, for instance, be used as a starting point in predicting transitions in boundary layers with computer simulations (computational fluid dynamics).
In the long-wave limit approximation, solutions are expressed in terms of generalized hypergeometric functions.

\end{abstract}


\end{frontmatter}


\section{Introduction}
\label{sec:1}
Stability and transition in shear flows and in boundary layers are mechanisms which need to be understood importantly for applications in mechanical and aerospace Engineering and in the atmospheric sciences \cite{KF,SH}. 
Both theoretical and experimental studies have been carried out by different researchers in order to improve our knowledge on the properties of mechanisms that govern transitions and instabilities in fluid flows and in boundary layers, and so, different types of instabilities that include linear inviscid, viscous and nonlinear instabilities, and transitions in fluid flows and in boundary layers were investigated and many more \cite{E,H,RS,SH} and references therein. 

With that goal in mind, analytical and numerical methods have been used to solve the Orr-Sommerfeld equation that governs the mechanisms of linear stabilities of fluid flows on finite domains (channels) , semi-infinite domains (e.g. boundary layers)  and infinite domains (e.g. wakes) \cite{E,H,SH}. The task becomes even more complicated when nonlinearities and turbulences are taken into consideration.  Due to the evolution of computers, this may be accomplished by performing computer simulations. However, obtaining initial input that allow computer simulations to converge to the correct solutions remains an important challenge.

For instance, Gregory et al. \cite{GS} performed numerical simulations and analyzed the linear inviscid stability of the three-dimensional boundary layer and applied it to the rotating disk flow but their results were biased. Brown \cite{B} extended the work by Gregory et al. \cite{GS} and included the viscous effects and applied the Orr-Sommerfeld equation to the rotating disc and swept-back wing. He used the temporal instability theory but his results did not match the observed values. Cebeci and Stewartson \cite{CS} used the spatial stability theory and solved the Orr-Sommerfeld equation on rotating disc profiles and their results were almost the same as Brown \cite{B}. 
Cooper and Carpenter \cite{CC1} investigated the stability of the rotating-disc boundary-layer flow over a compliant wall, and analyzed the the so called Type I and II instabilities. Reid \cite{R} derived an exact solution to the Orr-Sommerfeld equation for the plane Couette flow. Walker et al. \cite{WT} investigated a physically-based computational technique intended to estimate an initial guess for complex values of the wavenumber of a disturbance leading to the solution of Orr–Sommerfeld equation in a shear flow on a semi-infinite domain. However, the problem of stability and transition in shear flows and in boundary layers still is a challenging problem \cite{RS,SH}. 

Here, a different procedure which is rather analytical than numerical as in  Walker et al. \cite{WT} is proposed.  Two approaches are considered, the short-wave limit approximation and the long-wave limit approximation. As it will be shortly shown (section \ref{sec:3}), in the short-wave limit approximation, the wavelengths in the spanewise and streamwise directions are relatively short 
while in the long-wave limit approximation, on the other hand, the wavelengths in the spanewise and streamwise directions are relatively long.

In the short-wave approximation, the Orr-Sommerfeld equation is written as a system of two second order ordinary differential equations, and thus the eigenvalues are approximated using the WKB method and the eigenfunctions are approximated in terms of Green's functions. Their corresponding outer solutions are investigated since their properties (e.g. behavior at infinity) are easier to analyze than those of the solutions obtained by means of Green's functions.  For instance, instead of solving the Orr-Sommerfeld equation in the limit of infinite Reynolds numbers and using the fact that the mean flow velocity $\bar{U}(y) = 0$ in the boundary layer near the wall as described in Schmid and Henningson \cite{SH}, where critical layers and singularities are introduced in the equations while the Orr-Sommerfeld equation does not have any, we rather assume that $\bar{U}(y) = y, 0 <y < \infty$, approximate the eigenvalues using WKB methods and derive the approximate solutions in terms of Green's functions.

On one hand, the advantage of this procedure is that it works for any value of the Reynolds number, while on another hand, it can always be used to approximate eigenvalues and eigenfunctions of the Orr-Sommerfeld equation for any type of the mean flow velocity profile in three dimensions (3D) that can be approximated in two dimensions (2D) using the Squires' Theorem (transformation) as $\bar{U}(y) = ay^2+by+c$, where $a, b$ and $c$ are arbitrary constants. 
The approximated eigenvalues may be used as a starting point in predicting transitions in shear flow (e.g. boundary layers) in two or three dimensions with computer simulations (CFD).

In the long-wave limit, the Orr-Sommerfeld is reduced to a second order ordinary differential equation, and the solutions are written in terms of the generalized hypergeometric function. The obtained results can be usefull in geophysical fluid dynamics (GFD) where the size of the fluid disturbances can be of the same order as the radius of the earth or greater, for example, in the atmospheric boundary layer \cite{KF}.


\section{The three-dimensional linear stability model}\label{sec:2}

We consider an incompressible parallel flow in three dimensions with velocity components
\begin{equation}
u(x,y,z,t)=\bar{u}(y)+u^\prime (x,y,z,t),
\label{eq1}
\end{equation}
\begin{equation}
v(x,y,z,t)=v^\prime (x,y,z,t),
\label{eq2}
\end{equation}
\begin{equation}
w(x,y,z,t)=\bar{w}(y)+w^\prime (x,y,z,t),
\label{eq3}
\end{equation}
and pressure
\begin{equation}
p(x,y,z,t)=\bar{p}(y)+p^\prime (x,y,z,t),
\label{eq4}
\end{equation}
where the terms with a bar represent the mean flow quantities, while the terms with a prime represent small perturbation quantities or waves. We then substitute then (\ref{eq1}), (\ref{eq2}), (\ref{eq3}) and (\ref{eq4}) in the Navier-Stokes equation and neglect nonlinear terms (products of perturbation quantities). This yields the linearized  Navier-Stokes equations,
\begin{equation}
\frac{\partial u^\prime}{\partial t}+\bar{u}\frac{\partial u^\prime}{\partial x}+v\frac{d\bar{u}}{dy}+\bar{w}\frac{\partial u^\prime}{\partial z}=-\frac{\partial p^\prime}{\partial x}+\nu\nabla^2u^\prime+\nu\frac{d^2\bar{u}}{d y^2},
\label{eq5}
\end{equation}
\begin{equation}
\frac{\partial v^\prime}{\partial t}+\bar{u}\frac{\partial v^\prime}{\partial x}+\bar{w}\frac{\partial v^\prime}{\partial z}=-\frac{d\bar{p}}{dy}-\frac{\partial p^\prime}{\partial y}+\nu\nabla^2v^\prime,
\label{eq6}
\end{equation}
and 
\begin{equation}
\frac{\partial w^\prime}{\partial t}+\bar{u}\frac{\partial w^\prime}{\partial x}+v\frac{d\bar{w}}{dy}+\bar{w}\frac{\partial w^\prime}{\partial z}=-\frac{\partial p^\prime}{\partial z}+\nu\nabla^2w^\prime+\nu\frac{d^2\bar{w}}{d y^2},
\label{eq7}
\end{equation}
where $\nu$ is the fluid kinematic viscosity.

The full model describing the linear stability of the three-dimensional incompressible fluid flow comprises equations
(\ref{eq5}), (\ref{eq6}) and (\ref{eq7}) complemented by the continuity equation \cite{SH}

\begin{equation}
\frac{\partial u^\prime}{\partial x}+\frac{\partial v^\prime}{\partial y}+\frac{\partial w^\prime}{\partial z}=0.
\label{eq8}
\end{equation}
Applying the divergence operator to (\ref{eq5}), (\ref{eq6}) and (\ref{eq7}) gives

\begin{equation}
\nabla^2p^\prime +\frac{d^2 \bar{p}}{d y^2}=-2\frac{d \bar{u}}{d y}\frac{\partial v^\prime}{\partial x}-2\frac{d \bar{w}}{d y}\frac{\partial v^\prime}{\partial z}.
\label{eq9}
\end{equation}
We now eliminate the pressure term in (\ref{eq6}) by applying the Laplacian operator and combining with (\ref{eq9}). We then obtain

\begin{equation}
\nabla^2\frac{\partial v^\prime}{\partial t}+\bar{u}\nabla^2\frac{\partial v^\prime}{\partial x}+\bar{w}\nabla^2\frac{\partial v^\prime}{\partial z}-\frac{d \bar{u}}{d y}\frac{\partial v^\prime}{\partial x}-\frac{d \bar{w}}{d y}\frac{\partial v^\prime}{\partial z}-\nu\nabla^4v^\prime=0.
\label{eq10}
\end{equation}

We make all variables in (\ref{eq10}) non-dimensional  with respect to a typical reference speed $V$, a typical length scale $H$ in the $y$-direction and a typical length scale $L$ in both the streamwise direction ($x$-direction) and spanewise direction ($z$-direction). Thus, we find that
equation (\ref{eq10}) does not change but the Laplacian operator takes the form 

\begin{equation}
\nabla^2=r^2\frac{\partial^2}{\partial x^2}+\frac{\partial^2}{\partial y^2}+r^2\frac{\partial^2}{\partial z^2},
\label{eq11}
\end{equation}
where

\begin{equation}
r^2=\frac{H^2}{L^2}
\label{eq12}
\end{equation}
is the spatial aspect ratio, while on the other hand, the non-dimensional kinematic viscosity is given by 

\begin{equation}
\nu=\frac{1}{R}=\frac{L^2}{H^2}\frac{\nu^*}{LV}=\frac{1}{r^2}\frac{\nu^*}{LV},
\label{eq13}
\end{equation}
where $\nu^*$ is the dimensional kinematic viscosity and $R$ is the Reynolds number. And the Reynolds number $R$ can now be written in terms of the aspect ratio $r^2$ as

\begin{equation}
R=r^2\frac{LV}{\nu^*}.
\label{eq14}
\end{equation}

\section{The Orr-Sommerfeld equation: from three dimensions to two dimensions}
\label{sec:3}

We discuss the three dimensional model and reduce it to a two dimensional one. In stability theory of fluid flow, this procedure is known as Squire's theorem (or transformation) \cite{SH}.
We consider that the perturbations are small wavelike perturbations propagating in $xz$-plane with amplitude $\phi(y)$, and assume these perturbations are plane waves. Thus
\begin{equation}
v^\prime(x,y,z,t)=\phi(y) e^{i(\alpha x+\beta y-\omega t)},
\label{eq15}
\end{equation}
where $\alpha$ and $\beta$ are wavenumbers in streamwise and spanwise directions, and $\omega$ is the transient frequency of the waves. Substituting  (\ref{eq15}) in (\ref{eq10}) gives
\begin{equation}
\phi_{yyyy}-2r^2(\alpha^2+\beta^2)\phi_{yy}+r^4(\alpha^2+\beta^2)^2\phi-iR[(\alpha\bar{u}+\beta\bar{v}-\omega)(\phi_{yy}-r^2(\alpha^2+\beta^2)\phi)-(\alpha\bar{u}_{yy}+\beta\bar{v}_{yy})\phi]=0,
\label{eq16}
\end{equation}
where the subscript $y$ stands for differentiation with respect to $y$, and $\phi$, for instance, satisfies the boundary conditions
\begin{equation}
\phi(0)=\phi(\infty)=0\hspace{0.12cm}\mbox{and}\hspace{0.12cm} \phi_y(0)=\phi_y(\infty)=0.
\label{eq17}
\end{equation}
in the boundary layer. Or $\phi$ satisfies the boundary conditions
\begin{equation}
\phi(-\infty)=\phi(\infty)=0\hspace{0.12cm}\mbox{and}\hspace{0.12cm} \phi_y(-\infty)=\phi_y(\infty)=0.
\label{eq18}
\end{equation}
This boundary condition, for example, shall be important in investigating the stability of the two dimensional wake.

We importantly observe that (\ref{eq16}) is the famous Orr-Sommerfeld equation but slightly modified by the aspect ratio $r^2$.  Setting $\tilde{\alpha}=r\alpha$, $\tilde{\beta}=r\beta, \tilde{R}=rR$ and $\tilde{\omega}=r\omega$ gives
\begin{equation}
\phi_{yyyy}-2(\tilde{\alpha}^2+\tilde{\beta}^2)\phi_{yy}+(\tilde{\alpha}^2+\tilde{\beta}^2)^2\phi-i\tilde{R}[(\tilde{\alpha}\bar{u}+\tilde{\beta}\bar{v}-\tilde{\omega})(\phi_{yy}-(\tilde{\alpha}^2+\tilde{\beta}^2)\phi)-(\tilde{\alpha}\bar{u}_{yy}+\tilde{\beta}\bar{v}_{yy})\phi]=0,
\label{eq19}
\end{equation}
which is  the Orr-Sommerfeld equation.

Next, grouping together the terms according to $\phi_{yyyy}$, $\phi_{yy}$ and $\phi$ and rearranging terms, (\ref{eq16}) is 
written as
\begin{multline}
\phi_{yyyy}
-[2r^2(\alpha^2+\beta^2)+iR(\alpha\bar{u}+\beta\bar{v}-\omega)]\phi_{yy}\\+[r^4(\alpha^2+\beta^2)^2+ir^2R(\alpha\bar{u}+\beta\bar{v}-\omega)(\alpha^2+\beta^2)+iR(\alpha\bar{u}_{yy}+\beta\bar{v}_{yy})]\phi=0.
\label{eq20}
\end{multline}
We further let $k^2 =\alpha^2+\beta^2$ so that $\alpha = k \cos \theta$ and $\beta = k \sin \theta$  in (\ref{eq20}), where $\theta$ is the angle of orientation of the phase velocity in the $xz$-plane, then we obtain
\begin{multline}
\phi_{yyyy}-[2r^2k^2+iRk(\bar{u} \cos \theta+\bar{v} \sin \theta-\omega/k)]\phi_{yy}\\+[r^4k^4+ir^2Rk^3(\bar{u} \cos \theta+\bar{v} \sin\theta-\omega/k)+iRk(\bar{u}_{yy} \cos \theta+\bar{v}_{yy} \sin\theta)]\phi=0.
\label{eq21}
\end{multline}
Moreover setting $\bar{U}=\bar{u} \cos \theta+\bar{v} \sin \theta$ and $\bar{U}_{yy}=\bar{u}_{yy} \cos \theta+\bar{v}_{yy} \sin \theta$ in (\ref{eq19}) gives the two dimensional Orr-Sommerfeld equation

\begin{equation}
\phi_{yyyy}-[2r^2k^2+iRk(\bar{U}-\omega/k)]\phi_{yy}+[r^4k^4+ir^2Rk^3(\bar{U}-\omega/k)+iRk\bar{U}_{yy}]\phi=0.
\label{eq22}
\end{equation}

We consider two configurations as mentioned before. In one configuration, we assume that the spacial aspect ratio is large $r \gg1$. This configuration is representative for waves with relatively short wavelengths in the $xz$-plane since the wavelength in the streamwise direction $\lambda_x=2\pi/\tilde{\alpha}=2\pi/(r{\alpha})$ and that in the spanewise  $\lambda_z=2\pi/\tilde{\beta}=2\pi/(r{\beta})$ become shorter and shorter as the aspect ratio becomes larger and larger ($r\to\infty$). And so, this is a short-wave limit approximation. In that case, it is possible to write the two-dimensional Orr-Sommerfeld equation (\ref{eq20}) in terms of a system of two second ordinary differential equations. We will shortly see (section \ref{sec:4}) that WKB methods can be used to approximate eigenvalues.

In the other configuration, the aspect ratio is small, $0 < r \ll 1$ ($r\to0^+$). And so $\lambda_x$ and $\lambda_z$ become larger as the aspect ratio becomes small ($r\to 0^+$). Therefore, this configuration is representative for waves with long wavelength in the $xz$-plane. Hence, this a long-wave limit approximation, and we will illustrate it (section \ref{sec:5}) with two examples in which analytical solutions can be obtained.

\section{Short-wave limit approximation ($r\gg 1$)}\label{sec:4}

According to (\ref{eq14}),  the Reynolds number $R$ is proportional to $r^2$. Therefore, the last term in the coefficient for $\phi$ in (\ref{eq20}) which is proportional to $r^2$ is negligible compared to the other terms which are proportional to $r^4$. And so, it  can be dropped in the short-wave limit configuration where $r \gg1$.  
We then have
\begin{equation}
\phi_{yyyy}-[2r^2k^2+iRk(\bar{U}-\omega/k)]\phi_{yy}+[r^4k^4+ir^2Rk^3(\bar{U}-\omega/k)]\phi=0.
\label{eq23}
\end{equation}
We note that the assumption $r\gg 1$ is quite important (see section 4.1) since it helps us to write (\ref{eq21}) in a form allowing us to make use of WKB methods to approximate solutions for (\ref{eq20}). We further observe that (\ref{eq20}) will become  (\ref{eq21}) if $\bar{U}(y)$ is a linear function of $y$, e.g. Couette flow. 

Now let us consider the differential operator
\begin{equation}
\left[\frac{d^2}{dy^2}+Q(k,\omega,y)\right]\left[\frac{d^2}{dy^2}+P(k,\omega,y)\right],
\label{eq24}
\end{equation}
its expansion is 
\begin{multline}
\left[\frac{d^2}{dy^2}+Q(k,\omega,y)\right]\left[\frac{d^2}{dy^2}+P(k,\omega,y)\right]=\frac{d^4}{dy^4}+[P(k,\omega,y)+Q(k,\omega,y)]\frac{d^2}{dy^2}\\+2\frac{d P}{dy}(k,\omega,y)\frac{d}{dy}+\left[\frac{d^2 P}{dy^2}(k,\omega,y)+Q(k,\omega,y)P(k,\omega,y)\right].
\label{eq25}
\end{multline}
Applying this differential operator to $\phi$ gives 
\begin{equation}
\phi_{yyyy}+(P+Q)\phi_{yy}+2P_y\phi_y+(P_{yy}+QP)\phi=0.
\label{eq26}
\end{equation}
Comparing this with (\ref{eq21}) yields
\begin{equation}
P+Q=-[2r^2k^2+iRk(\bar{U}-\omega/k)],
\label{eq27}
\end{equation}
\begin{equation}
QP=r^4k^4+ir^2Rk^3(\bar{U}-\omega/k)
\label{eq28}
\end{equation}
and 
\begin{equation}
P=\mbox{constant}.
\label{eq29}
\end{equation}
This gives
\begin{equation}
P=-r^2k^2.
\label{eq30}
\end{equation}
and
\begin{equation}
Q=-[r^2k^2+iRk(\bar{U}-\omega/k)]=-r^2[k^2+i\chi k(\bar{U}-\omega/k)],
\label{eq31}
\end{equation}
where according to (\ref{eq14}), $\chi=R/r^2={LV}/{\nu^*}$.
Hence, setting in (\ref{eq24})
\begin{equation}
\phi_{yy}+P\phi=\Psi.
\label{eq32}
\end{equation}
implies that $\Psi$ has to satisfy
\begin{equation}
\Psi_{yy}+Q\Psi=0.
\label{eq33}
\end{equation}

We use the following theorem to establish the boundary conditions for (\ref{eq33}).
\begin{thm}
\begin{enumerate} 
\item Consider equation (\ref{eq22}) with boundary conditions $\phi(a) =\phi(\infty)  = 0$  and  $\phi_y(a)=\phi_y(\infty) = 0$, $ a\ge0$.
\begin{enumerate}
\item If $\phi$ and $\Psi$ are solutions for (\ref{eq32}) and (\ref{eq33}) respectively, then $\Psi(\infty) = 0$.
\item There exists some constants $\delta_1>0$ and $\delta_2$ such that $\Psi(a+\delta_1)=\delta_2$.
\item And if further the aspect ratio $r\gg1$ then $\Psi(a + \delta_1) \le\delta_2$.
\end{enumerate}
\item Consider equation (\ref{eq21}) with boundary conditions $\phi(-\infty) =\phi(a)  = 0$  and  $\phi_y(-\infty)=\phi_y(a) = 0$, $ a\le0$.
\begin{enumerate}
\item If $\phi$ and $\Psi$ are solutions for (\ref{eq32}) and (\ref{eq33}) respectively, then $\Psi(-\infty) = 0$.
\item There exists some constants $\delta_1>0$ and $\delta_2$ such that $\Psi(a-\delta_1)=\delta_2$.
\item And if further the aspect ratio $r\gg1$ then $\Psi(a - \delta_1) \le\delta_2$.
\end{enumerate}
\end{enumerate}
\label{thm1}
\end{thm}

\begin{pf}
\begin{enumerate} 
\item 
\begin{enumerate}
\item We first observe that the solution $\phi$ for (\ref{eq32}) is $\phi(y) =\phi_h(y) +\phi_p(y)= A e^{rky} + B e^{-rky} + \phi_p(y)$, where $A$ and $B$ are constants and $\phi_h$ and $\phi_p$ are the homogeneous solution and particular solution respectively. Applying the method of undetermined coefficients to (\ref{eq32}), the particular solution has to take the form $\phi_p(y) = f(y)\Psi(y)$ where $f(y)$ is a function that has to be chosen in order to make $\phi_h$ and $\phi_p$ independent. But the
form of $P(y)$ and that of $Q(y)$ given by (\ref{eq29}) and (\ref{eq30}) respectively, indicate that $\phi_h$ and $\phi_p$ will always be independent for all $y$. And so $\phi_p=\varepsilon\Psi(y)$ where $\varepsilon$ is some constant. Moreover $A$ must vanish in order $\phi(y)$ to be finite as $y\to\infty$. And so $\phi(y) = B e^{-rky} + \varepsilon \Psi(y)$. Therefore, $\phi(\infty) = \lim_{y\to\infty}(B e^{-rky} + \eta\Psi(y)) = \varepsilon \lim_{y\to\infty}\Psi(y) = 0$. Hence, $\Psi(\infty) = 0$ since $\varepsilon$ is a constant.
\item $\phi(a+\delta_1)=\lim_{y\to a+\delta_1}( B e^{-rky} + \varepsilon \Psi(y))= B e^{-rk(a+\delta_1)} + \varepsilon \Psi(a+\delta_1)=\Delta $, some constant. And so, $\Psi(a+\delta_1)=\frac{\Delta}{\varepsilon}-\frac{B e^{-rk(a+\delta_1)}}{\varepsilon}=\delta_2$.
\item If further the aspect ratio $r\gg1$, then $\frac{B e^{-rk(a+\delta_1)}}{\varepsilon}\to0$. Hence, $\Psi(a + \delta_1) \le\delta_2$.
\end{enumerate}
\item
\begin{enumerate}
\item From 1. (a), $\phi(y) =\phi_h(y) +\phi_p(y)= A e^{rky} + B e^{-rky} + \varepsilon\Psi(y)$, where $A$, $B$ and $\varepsilon$ are constants. And so $B$ must vanish in order $\phi(y)$ to be finite as $y\to -\infty$. This gives $\phi(y) = A e^{rky} +\varepsilon\Psi(y)$. Therefore, $\phi(-\infty) = \lim_{y\to -\infty}(B e^{rky} +\varepsilon\Psi(y)) = \varepsilon \lim_{y\to -\infty}\Psi(y) = 0$. Hence, $\Psi(-\infty) = 0$ since $\varepsilon$ is a constant.
\item $\phi(a-\delta_1)=\lim_{y\to a-\delta_1}( B e^{rky} + \varepsilon \Psi(y))= B e^{rk(a-\delta_1)} + \varepsilon \Psi(a-\delta_1)=\Delta $, some constant. And so, $\Psi(a-\delta_1)=\frac{\Delta}{\varepsilon}-\frac{B e^{rk(a-\delta_1)}}{\varepsilon}=\delta_2$.
\item If further the aspect ratio $r\gg1$, then $\frac{B e^{rk(a-\delta_1)}}{\varepsilon}\to0$. Hence, $\Psi(a -\delta_1) \le\delta_2$.
\end{enumerate}
\end{enumerate}
\end{pf} 

\begin{clr}
Consider equation (\ref{eq21}) with boundary conditions $\phi(-\infty) = \phi(\infty) = 0$ and $\phi_y(-\infty) = \phi_y(\infty) = 0$. If $\phi$ and $\Psi$ are solutions for (\ref{eq32}) and (\ref{eq33}) respectively, then $\Psi(-\infty) = \Psi(\infty) = 0$.
\label{clr1}
\end{clr}

\begin{pf}
Corollary \ref{clr1}'s proof follows from Theorem 1. To prove Corollary \ref{clr1}, we set $a=0$, $A=B$ and let $\delta_1\to0$ in Theorem 1 so that we obtain $\Psi(-\infty) = \Psi(\infty) = 0$.
\end{pf}

\subsection{WKB approximation for eigenvalues on a semi-infinite domain in the short wave-limit approximation}\label{subsec:4.1}
Using Theorem \ref{thm1}, in the boundary layer, we solve
\begin{equation}
\phi_{yy}+P\phi=\Psi.
\label{eq34}
\end{equation}
and
\begin{equation}
\Psi_{yy}+Q\Psi=0,
\label{eq35}
\end{equation}
with boundary conditions
\begin{equation}
\phi_y(0)=\phi_y(\infty) = 0
\label{eq36}
\end{equation}
and
\begin{equation}
\Psi(0+\delta_1)=\delta_2 \hspace{.12cm} \mbox{and}\hspace{.12cm} \Psi(\infty) = 0,
\label{eq37}
\end{equation}
where $\delta_1$ is a very small positive constant ($\delta_1\to0^+$), and as before, $P = -r^2k^2$ and $Q = -r^2[k^2+i\chi k(\bar{U}(y)-\omega/k)]$. 

Following Theorem \ref{thm1}, $\delta_2$ may be set to zero. In that case, we can restrict the domain to $R^+ = (0,\infty)$, and thus use the WKB method described in \ref{AppB} to approximate the eigenvalues $k$. We write $Q$ as,
\begin{equation}
Q(k,\omega,y) = -r^2[k^2+i\chi k(\bar{U}(y)-\omega/k)]=i r^2 \chi k[\lambda(k,\omega)-\bar{U}(y)],
\label{eq38}
\end{equation}
where $\lambda(k,\omega)=(i\chi\omega-k^2)/i\chi k$. Now, (\ref{eq35}) can be rewritten as the Schrodinger
equation \cite{BO},
\begin{equation}
\epsilon^2\Psi_{yy}-[V-E]\Psi=0,
\label{eq39}
\end{equation}
where $E=i\chi k\lambda(k,\omega), V(y)=i\chi k\bar{U}(y)$ and $\epsilon^2=1/r^2\to0$ is a small constant since $r^2\gg1$.

It is shown in the \ref{AppB} that using the WKB method, the eigenvalues do satisfy
\begin{equation}
\int\limits_0^b\sqrt{E-V(y)}dy\sim\left(n-\frac{1}{4}\right)\frac{\pi}{r}, r\gg1, n = 0, 1, 2, \cdots,
\label{eq40}
\end{equation}
where the limits of integration $y=0$ and $y=b$ are the turning points of $V-E$. On substituting $E=i\chi k\lambda(k,\omega)$ and $V(y)=i\chi k\bar{U}(y)$ in (\ref{eq40}) gives
\begin{equation}
\int\limits_0^b\sqrt{\lambda(k,\omega)-\bar{U}(y)}dy\sim\frac{1}{\sqrt{i\chi k}}\left(n-\frac{1}{4}\right)\frac{\pi}{r}, r\gg1, n = 0, 1, 2, \cdots.
\label{eq41}
\end{equation}

We consider two configurations and approximate the eigenvalues using WKB method \cite{BO}. In the first configuration (configuration 1), the background mean flow is given by $\bar{U}(y)=by+c, 0<y<\infty$, where $b$ and $c$ are constants, while in the second configuration (configuration 2), the background mean flow is given by $\bar{U}(y)=ax^2+by+c, 0<y<\infty$,  where $a,b$ and $c$ are constants.

\subsubsection{Configuation 1 ($r\gg1$): The linear background flow, $\bar{U}(y)=by+c, 0<y<\infty$}
\label{subsub:4.1.1}
Let us consider the linear background flow $\bar{U}(y)=by+c,0 < y <\infty$, where $b$ and $c$ are constants.
The turning points are $y=0$ and $y_= (\lambda(k,\omega)-c)/b$. Applying (\ref{eq41}) gives the dispersion relation
\begin{equation}
\frac{2}{3}\left(\frac{i\chi\omega-k^2}{i\chi k}-c\right)^{3/2}\sim\frac{b}{\sqrt{i\chi k}}\left(n-\frac{1}{4}\right)\frac{\pi}{r}.
\label{eq42}
\end{equation}
Hence,
\begin{equation}
\omega_n\sim-i\frac{k_n^2}{\chi}+ck_n+\left(\frac{3bk_n}{2\sqrt{i\chi}}\right)^{2/3}\left(n-\frac{1}{4}\right)^{2/3}\left(\frac{\pi}{r}\right)^{2/3}.
\label{eq43}
\end{equation}

If the flow is steady (e.g. laminar boundary layer), $\omega= 0$ and the eigenvalues $k_n$ satisfy 
\begin{equation}
i\frac{k_n^2}{\chi}-ck_n-\left(\frac{3bk_n}{2\sqrt{i\chi}}\right)^{2/3}\left(n-\frac{1}{4}\right)^{2/3}\left(\frac{\pi}{r}\right)^{2/3}=0
\label{eq44}
\end{equation}
which can be solved using basic numerical methods. In the special case where $\bar{U}(y)=y$ ($b=1$ and $c=0$), (\ref{eq44}) can be explicitly solved to obtain
\begin{equation}
k_n\sim\mp(1+i)\frac{\sqrt{2}}{2}\left(\frac{3\chi}{2r}\right)^{1/2}\left(n-\frac{1}{4}\right)^{1/2}\pi^{1/2}
\label{eq45}
\end{equation}
If for example $\theta=\pi/6$, then the streamwise wavenumber is
\begin{equation}
\alpha_n=k_n\cos\theta\sim\mp(1+i)\frac{\sqrt{2}}{4}\left(\frac{3\chi}{2r}\right)^{1/2}\left(n-\frac{1}{4}\right)^{1/2}\pi^{1/2},
\label{eq46}
\end{equation}
while the spanewise wavenumber is 
\begin{equation}
\beta_n=k_n\sin\theta\sim\mp(1+i)\frac{\sqrt{6}}{4}\left(\frac{3\chi}{2r}\right)^{1/2}\left(n-\frac{1}{4}\right)^{1/2}\pi^{1/2}.
\label{eq47}
\end{equation}
\subsubsection{Configuration 2  ($r\gg1$): The quadratic background flow $\bar{U}(y)=ay^2+by+c, 0<y<\infty$}
\label{subsub:4.1.3}

Here, we consider the quadratic background flow $\bar{U}(y)=ay^2+by+c, 0<y<\infty$, where $a,b$ and $c$ are constants. We rewrite (\ref{eq41}) as 
\begin{equation}
\int\limits_{-b/(2a)}^{\sqrt{\Lambda}-b/(2a)}\sqrt{\Lambda-\mu^2}d\mu\sim\frac{1}{\sqrt{ia\chi k}}\left(n-\frac{1}{4}\right)\frac{\pi}{r},
\label{eq48}
\end{equation}
where $\Lambda=[\lambda+b^2/(4a)-c]/a$ and $\mu=y+(b/2a)$.
This gives the dispersion relation
\begin{equation}
\frac{1}{4}\left(\frac{i\chi\omega-k^2}{i\chi k}+\frac{b^2}{4a}-c\right)\sim\sqrt{\frac{a}{i\chi k}}\left(n-\frac{1}{4}\right)\frac{1}{r}.
\label{eq49}
\end{equation}
Therefore,  we obtain the asymptotic approximation for $\omega$,
\begin{equation}
\omega_n\sim-i\frac{k_n^2}{\chi}-k_n\left(\frac{b^2}{4a}-c\right)+\frac{4}{r}\sqrt{\frac{ak_n}{i\chi}}\left(n-\frac{1}{4}\right).
\label{eq50}
\end{equation}

In the steady flow where $\omega= 0$, $k_n$ satisfies 
\begin{equation}
i\frac{k_n^2}{\chi}+k_n\left(\frac{b^2}{4a}-c\right)-\frac{4}{r}\sqrt{\frac{ak_n}{i\chi}}\left(n-\frac{1}{4}\right)=0.
\label{eq51}
\end{equation}
which can be solved using basic numerical methods. In the special case where $\bar{U}(y)=y^2$ ($a=1,b=0$ and $c=0$)
\begin{equation}
k_n\sim(\sqrt{3}+i)\left(\frac{2\chi}{r^2}\right)^{1/3}\left(n-\frac{1}{4}\right)^{2/3}.
\label{eq52}
\end{equation}
If $\theta=-\pi/4$ for example,  the streamwise and spanewise wavenumbers are respectively given by
\begin{equation}
\alpha_n=k_n\cos\theta\sim(\sqrt{3}+i)\frac{\sqrt{2}}{2}\left(\frac{2\chi}{r^2}\right)^{1/3}\left(n-\frac{1}{4}\right)^{2/3}
\label{eq53}
\end{equation}
while
\begin{equation}
\beta_n=k_n\sin\theta\sim-(\sqrt{3}+i)\frac{\sqrt{2}}{2}\left(\frac{2\chi}{r^2}\right)^{1/3}\left(n-\frac{1}{4}\right)^{2/3}.
\label{eq54}
\end{equation}

\subsection{Green's function approximation for eigenfunctions}\label{subsec:4.2}

In this section, we approximate the eigenfunctions $\phi_n$ by means of Green's functions. We first solve (\ref{eq35}) for $\Psi$ and then use Green's functions to solve the inhomogeneous equation (\ref{eq34}).  In that case, solutions to (\ref{eq34})  will be given by 
\begin{equation}
\phi_n(y)=\int\limits_0^\infty G(y,\xi)\Psi_n(\xi)d\xi,
\label{eq55}
\end{equation}
where $G(y,\xi)$ are Green's functions. Green's functions associated with the boundary value problem defined by (\ref{eq34}) and (\ref{eq36}) are given by (\ref{eqA8}) (see \ref{AppA}).

\subsubsection{Configuration 1  ($r\gg1$): $\bar{U}(y)=by+c,0<y<\infty$}\label{subsub:4.2.1}

We consider that the background flow is linear, $\bar{U}=by+c$, where $b$ and $c$ are constants as before. We first solve, for  $\Psi$,  
\begin{equation}
\Psi_{yy}- i r^2\chi k [by+c-\lambda(k,\omega)]\Psi=0
\label{eq56}
\end{equation}
subject to boundary conditions
\begin{equation}
\Psi(0+\delta_1)=0, \delta_1 \to 0^{+} \hspace{.15cm} \mbox{and} \hspace{.15cm} \Psi(\infty)=0,
\label{eq57}
\end{equation}
as before. We make the change of variable $\eta= (ir^2\chi kb)^{1/3} (y-\tilde{\lambda})$, where $\tilde{\lambda}=-ir^2\chi k[c-\lambda(k,\omega)]/b$. Rearranging terms gives the Airy equation \cite{AB}
\begin{equation}
\Psi_{\eta\eta}-\eta\Psi=0
\label{eq58}
\end{equation}
which has solution 
\begin{equation}
\Psi(\eta)=D_1\mbox{Ai}(\eta)+D_2\mbox{Bi}(\eta),
\label{eq59}
\end{equation}
where $D_1$ and $D_2$ are constants \cite{AB}. 

A solution satisfying the boundary conditions (\ref{eq57}) is thus given by 
\begin{equation}
\Psi_n(y)=\mbox{Ai}[ (ir^2\chi k_n b)^{1/3} (y-\tilde{\lambda}_n)],  n=0,1,2,\cdots.
\label{eq60}
\end{equation}
Hence, applying (\ref{eq55}) gives
\begin{multline}
\phi_n(y)\sim\frac{e^{-rk_n y}}{rk_n}\int\limits_0^y \cosh(rk_n \xi)\mbox{Ai}[ (ir^2\chi k_n b)^{1/3} (\xi-\tilde{\lambda}_n)]d\xi\\+
\frac{\cosh(rk_n y)}{rk_n}\int\limits_y^\infty e^{-rk_n\xi}\mbox{Ai}[ (ir^2\chi k_n b)^{1/3} (\xi-\tilde{\lambda}_n)]d\xi, n=0,1,2,\cdots.
\label{eq61}
\end{multline}

\subsubsection{Configuration 2  ($r\gg1$): $\bar{U}(y)=ay^2+by+c, 0<y<\infty$}\label{subsub:4.2.2}
The background mean flow is $\bar{U}(y)=ay^2+by+c, 0<y<\infty$, where $a,b$ and $c$ are constants as before. Then, $\Psi$ satisfies 
\begin{equation}
\Psi_{yy}- i r^2\chi k [ay^2+by+c-\lambda(k,\omega)]\Psi=0
\label{eq62}
\end{equation}
subject to boundary conditions
\begin{equation}
\Psi(0+\delta_1)=0, \delta_1 \to 0^{+} \hspace{.15cm} \mbox{and} \hspace{.15cm} \Psi(\infty)=0
\label{eq63}
\end{equation}
as before. Now, setting $\tilde{\lambda}=r\sqrt{i\chi k}[\lambda+b^2/(4a)-c]/\sqrt{a}$, making the change of variable $\eta= (ir^2\chi ka)^{1/4} [y+(b/2a)]$, and rearranging terms gives the Hermite equation \cite{AB}
\begin{equation}
\Psi_{\eta\eta}(\eta)-[\eta^2-\tilde{\lambda}]\Psi(\eta)=0
\label{eq64}
\end{equation}
which has eigen-solutions 
\begin{equation}
\Psi_n(\eta)=\frac{\pi^{-1/4}e^{-\eta^2}}{\sqrt{2^n 2!}}H_n(\eta), n=0,1,2,\cdots,
\label{eq65}
\end{equation}
where $H_n$ are Hermite polynomials of $n$ \cite{AB}. 

The boundary condition (\ref{eq63}) shall be satisfied if and only if $n$ is odd. In that case, solutions to  (\ref{eq62})- (\ref{eq63}) take the form
\begin{equation}
\Psi_{m}(y)=\frac{\exp\left(-(ir^2\chi k_{_{2m+1}})^{\frac{1}{2}} y^2\right)}{\pi^{\frac{1}{4}}\sqrt{2^{2m+1} 2!}}H_{2m+1}\left((ir^2\chi k_{2m+1}a)^{\frac{1}{4}} \left (y+\frac{b}{2a}\right)\right), m=0,1,2,\cdots,
\label{eq66}
\end{equation}
where $H_{2m+1}$ are Hermite polynomials of $2m+1$. To obain the eigen-solutions $\phi_m(y)$, we apply (\ref{eq55}) as as in configuration 1.

\subsection{Outer solution approximation for the eigenfunctions}\label{subsec:4.3}

The analysis of the eigenfunctions obtained by means of Green's functions does not seem to be an easy task due to the complexities involved in the computations of the integrals. In order to understand the behaviors of these solutions we can instead look at the outer solutions. Outer solutions are valid when $\epsilon=1/r\to 0$ is a small parameter and give insight into the behavior of the solution for large argument $y\gg1$ (Bender and Orszag \cite{BO}). In the short-wave limit approximation, the homogeneous solution is proportional to $e^{-rky} = e^{-ky/\epsilon}$. Therefore, we expect the homogeneous solution to quickly vanish with $y$ as $\epsilon\to0$. In that case, the outer solution would be accurate even for small values of $y$ of order $O(\epsilon/k),\epsilon\to0$. Hence, outer solutions shall be valid on $\mathbb{R}^+$ as long as $\epsilon\to0$.

For small enough $\epsilon$, we obtain from (\ref{eq32}) that $P\phi\sim\Psi$. This gives
\begin{equation}
\phi(y)\sim \frac{\Psi(y)}{P}=-\frac{\Psi(y)}{r^2k^2}=-\frac{\epsilon^2\Psi(y)}{k^2}.
\label{eq67}
\end{equation}
This means that the solution $\phi$ of the Orr-Sommerfeld equation (\ref{eq16}) is driven by $\Psi$ whenever $\epsilon\to0$ or the aspect ratio is large ($r\gg1$) since the homogeneous solution which is proportional to $e^{-rky} = e^{ky/\epsilon}$ rapidly vanishes with $y$. Thus, the constant $\varepsilon$ in Theorem \ref{thm1} should be $\varepsilon=r^2/k^2$.  We look at two special cases of configurations 1 and 2 (see section \ref{subsec:4.2}), for which the outer solutions and the eigenvalues can explicitly be obtained.

\subsubsection{Case 1: $\bar{U}(y)=y, 0<y<\infty$}

This case corresponds to the configuration 1 in section \ref{subsec:4.2} with the constants $b=1$ and $c=0$. We then have $\tilde{\lambda}=\lambda$ in (\ref{eq60}). We also consider that $\omega=0$ (laminar flow) and let the phase velocity angle orientation $\theta=0$ so that $k=\alpha$ the streamwise wavenumber. Using (\ref{eq60}) and (\ref{eq67}) gives the outer solution 
\begin{align}
\phi_n(y)&\sim-\frac{\Psi_n(y)}{r^2\alpha_n^2}=-\frac{\mbox{Ai}\{(ir^2\chi \alpha_n)^{1/3} [y-\lambda_n(\alpha_n)]\}}{r^2\alpha_n^2}
\nonumber \\ &=-\frac{\epsilon^2\Psi_n(y)}{\alpha_n^2}=-\frac{\epsilon^2}{\alpha_n^2}\mbox{Ai}\{(iR\alpha_n)^{1/3} [y-\lambda_n(\alpha_n)]\},  n=0,1,2,\cdots, 
\label{eq68}
\end{align}
where $\lambda_n$ (see section \ref{subsec:4.1}) is given by 
\begin{equation}
\lambda_n(\alpha_n)=\left[\frac{3}{2r}\frac{1}{\sqrt{i\chi\alpha_n}}\left(n-\frac{1}{4}\right)\pi\right]^{2/3}=\left[\frac{3}{2}\frac{1}{\sqrt{iR\alpha_n}}\left(n-\frac{1}{4}\right)\pi\right]^{2/3}, 
\label{eq69}
\end{equation}
with, see equation (\ref{eq45}),
\begin{multline}
\alpha_n\sim\mp(1+i)\frac{\sqrt{2}}{2}\left(\frac{3\chi}{2r}\right)^{1/2}\left(n-\frac{1}{4}\right)^{1/2}\pi^{1/2}\\=\mp\epsilon(1+i)\frac{\sqrt{2}}{2}\left(\frac{3R\epsilon}{2}\right)^{1/2}\left(n-\frac{1}{4}\right)^{1/2}\pi^{1/2}.
\label{eq70}
\end{multline} 

Some plots of the outer eigen-solution $\phi_n$ as a function of $y$ are shown on Figure \ref{fig1}. It is seen that the amplitude of $\phi_n$ increases as the Reynolds number $R$ decreases.

\subsubsection{Case 2: $\bar{U}(y)=y^2, 0<y<\infty$}

This case corresponds to the configuration 2 in section \ref{subsec:4.2} with the constants $a=1, b=0$ and $c=0$ in (\ref{eq69}). We consider that $\omega=0$ (laminar flow) and let $\theta=0$ so that $k=\alpha$ the streamwise wavenumber. We use (\ref{eq69}) and (\ref{eq70}), and then we obtain 
\begin{align}
\phi_{m}(y)&\sim-\frac{\Psi_{m}(y)}{r^2\alpha_{m}^2}=-\frac{\exp\left(-(ir^2\chi \alpha_{_{2m+1}})^{\frac{1}{2}} y^2\right)}{\pi^{\frac{1}{4}}r^2\alpha_{2m+1}^2\sqrt{2^{2m+1} 2!}}H_{2m+1}\left((ir^2\chi \alpha_{2m+1})^{\frac{1}{4}} y\right) \nonumber \\&=-\frac{\epsilon^2\Psi_{m}(y)}{\alpha_m^2}=-\frac{\epsilon^2}{\alpha_{m}^2}\frac{\exp\left(-(iR \alpha_{_{2m+1}})^{\frac{1}{2}} y^2\right)}{\pi^{\frac{1}{4}}\sqrt{2^{2m+1} 2!}}H_{2m+1}\left((iR \alpha_{2m+1})^{\frac{1}{4}} y\right), m=0,1,2,\cdots,
\label{eq71}
\end{align}
where, as before, $H_{2m+1}$ are Hermite polynomials of order $2m+1$ \cite{AB}, and,
\begin{equation}
\alpha_{m}=k_m\sim(\sqrt{3}+i)\left(\frac{2\chi}{r^2}\right)^{1/3}\left(2m+\frac{3}{4}\right)^{2/3}= (\sqrt{3}+i)\left(2R\epsilon^4\right)^{1/3}\left(2m+\frac{3}{4}\right)^{2/3}.
\label{eq72}
\end{equation}

Some results are shown in Figure \ref{fig2}. It is seen that the amplitude of $\phi_n$ increases as the Reynolds number $R$ decreases as in case 1 where $\bar{U}(y)=y$.

\subsection{Stability of the two-dimensional wake in the short-wave limit}\label{subsec:4.4}

We write our solutions in terms of the hypergeometric function $_2F_1$. For reference, we shall first define the generalized hypergeometric function. Functions of this type are also used in section \ref{sec:5}.

\begin{dfn}The generalized hypergeometric function, denoted as $_pF_q$, is a special function given by the series \cite{AB,NI}
\begin{equation}
_p F_q(a_1, a_2,\cdots,a_p;b_1, b_2, \cdots, b_q; x)=\sum\limits_{n=0}^{\infty}\frac{(a_1)_n (a_2)_n\cdots (a_p)_n}{(b_1)_n (b_2)_n\cdots (b_q)_n}\frac{x^n}{n!},
\label{eq73}
\end{equation} 
where $a_1, a_2,\cdots,a_p$ and $;b_1, b_2, \cdots, b_q$ are arbitrary constants, $(\vartheta)_n=\Gamma(\vartheta+n)/\Gamma(\vartheta)$ (Pochhammer's notation \cite{AB}) for any complex $\vartheta$, with $(\vartheta)_0=1$, and $\Gamma$ is the standard gamma function \cite{AB,NI}.
\label{dfn1}
\end{dfn}

In this section, we apply our method to a problem where $\phi(y)$ is defined on infinite domain rather than on a semi-infinite domain as in section \ref{subsec:4.3}. We investigate the stability of a two-dimensional wake. We let the background flow resemble a Bell shape function given by

\begin{equation}
\bar{U}(y)=U_0\ \mbox{sech}^2(\varphi y), -\infty<y<+\infty,
\label{eq74}
\end{equation}
where $U_0$ and $\varphi$ are positive constants. 

We make use of Corollary \ref{clr1} to establish the values of $\Psi$ at the boundaries of the domain. Hence, we solve
\begin{equation}
\phi_{yy}-r^2k^2\phi=\Psi.
\label{eq75}
\end{equation}
and
\begin{equation}
\Psi_{yy}-ir^2\chi k[U_0\ \mbox{sech}^2(\varphi y)-\lambda(k,\omega)]\Psi=0,
\label{eq76}
\end{equation}
subject to boundary conditions
\begin{equation}
\phi_y(-\infty)=\phi_y(\infty) = 0
\label{eq77}
\end{equation}
and
\begin{equation}
\Psi(-\infty)=\Psi(\infty) = 0,
\label{eq78}
\end{equation}
where as before $\chi=R/r^2$ and $\lambda=(i\chi\omega-k^2)/(i\chi k)$.

The Sturm-Liouvile problem defined by the Schrodinger equation (\ref{eq76}) with the boundary condition (\ref{eq78}) is a well known problem in quantum mechanics, see for example the Appendix A in Nyengeri \cite{NY} for details about the derivation of the solution. $\Psi_n$ is thus given by
\begin{equation}
\Psi_{n}(y)=\frac{_2 F_1\left[\frac{1}{2}+\vartheta_n+\sqrt{\frac{-ir^2k_n\chi U_0}{\varphi^2}+\frac{1}{4}},  \frac{1}{2}+\vartheta_n-\sqrt{\frac{-ir^2k_n\chi U_0}{\varphi^2}+\frac{1}{4}}; 1+\vartheta_n; \sigma(y)\right]}{\cosh^{\vartheta_n}(\varphi y)}, n=0,1,2,\cdots,
\label{eq79}
\end{equation}
where $ \sigma(y)=1-\tanh(\varphi y)$ and $ \vartheta_n=-(n+1/2)+[{-ir^2k_n\chi U_0}/{\varphi^2}+{1}/{4}]^{1/2}$.

We also have 
\begin{equation}
 \lambda(\omega_n, k_n)=-\varphi^2\left[\left(n+\frac{1}{2}\right)-\sqrt{\frac{-ir^2k_n\chi U_0}{\varphi^2}+\frac{1}{4}}\right].
\label{eq80}
\end{equation}
Using the fact that $\lambda=(i\chi\omega-k^2)/(i\chi k)$, we obtain the dispersion relation 

\begin{equation}
\omega_n\sim-i\frac{k_n^2}{\chi}-k_n\varphi^2\left[\left(n+\frac{1}{2}\right)-\sqrt{\frac{-ir^2k_n\chi U_0}{\varphi^2}+\frac{1}{4}}\right].
\label{eq81}
\end{equation}

Green's functions associated with the boundary value problem defined by (\ref{eq75}) and (\ref{eq77}) are given by (\ref{eqA10}) in \ref{AppA}.
Thus, applying (\ref{eq55}) gives the asymptotic eigenfunctions
\begin{multline}
\phi_{n}(y)=\frac{e^{rky}}{2rk}\int\limits_{-\infty}^{y}e^{-rk\xi}\frac{_2 F_1\left[\frac{1}{2}+\vartheta_n+\sqrt{\frac{-ir^2k_n\chi U_0}{\varphi^2}+\frac{1}{4}},  \frac{1}{2}+\vartheta_n-\sqrt{\frac{-ir^2k_n\chi U_0}{\varphi^2}+\frac{1}{4}}; 1+\vartheta_n; \sigma(\xi)\right]}{\cosh^{\vartheta_n}(\varphi \xi)}d\xi\\\hspace{1.3cm}+\frac{e^{-rky}}{2rk}\int\limits_{y}^{+\infty}e^{rk\xi}\frac{_2 F_1\left[\frac{1}{2}+\vartheta_n+\sqrt{\frac{-ir^2k_n\chi U_0}{\varphi^2}+\frac{1}{4}},  \frac{1}{2}+\vartheta_n-\sqrt{\frac{-ir^2k_n\chi U_0}{\varphi^2}+\frac{1}{4}}; 1+\vartheta_n; \sigma(\xi)\right]}{\cosh^{\vartheta_n}(\varphi \xi)}d\xi,\\ n=0,1,2,\cdots.
\label{eq82}
\end{multline}

\section{Long-wave limit approximation ($r\to0^+$) on a semi-infinite domain}\label{sec:5}

In this section, we consider the  semi-infinite domain $0<y<\infty$. As seen in section \ref{subsec:4.3}, the analysis of the solutions obtained by means of Green's functions is not an easy task, here we use a different procedure. However, in the long-wave limit approximation, WKB methods cannot be applied, but (\ref{eq20}) can be reduced to a form that allows us to readily obtain solutions in terms of hypergeometric functions whose properties are known. 

In the long-wave limit approximation ($r\to0^+$),  the terms proportional to $r^4$ in the coefficient of $\phi$ in (\ref{eq20}) are quite negligible in amplitude compared to those proportional to $r^2$, and so they can be dropped. Therefore, (\ref{eq20}) reduces to

\begin{equation}
\phi_{yyyy}-[2r^2k^2+ir^2\chi k(\bar{U}-\omega/k)]\phi_{yy}+ir^2\chi k\bar{U}_{yy}\phi=0.
\label{eq83}
\end{equation}
with boundary conditions
\begin{equation}
\phi(0)=\phi(\infty)=0\hspace{0.12cm}\mbox{and}\hspace{0.12cm} \phi_y(0)=\phi_y(\infty)=0.
\label{eq84}
\end{equation}
in the boundary layer. 

Here, we consider two velocity mean profiles, the linear velocity mean profile $\bar{U}(y) = by + c$ and the quadratic velocity mean profile $\bar{U}(y) =\delta y^2 + by + c$ and solve the boundary value problem (\ref{eq91}) - (\ref{eq92}). For the quadratic mean flow profile, we
assume that $\delta$ is a small constant. This implies that $\bar{U}_{yy}=\delta$ is also small, and consequently the third term involving $\bar{U}_{yy}=\delta$ may be dropped.

\subsection{Configuration 3 ($r\to0^+$):  $\bar{U}(y) = by + c, 0<y<\infty$}
We have $\bar{U}_{yy}=0$. Then equation (\ref{eq83}) becomes
\begin{equation}
\phi_{yyyy}-[2r^2k^2+ir^2\chi k(by+c-\omega/k)]\phi_{yy}=0.
\label{eq85}
\end{equation}
Setting $\Psi=\phi_{yy}, \lambda=[\omega/k-c-2k/(i\chi)]/b$ and making the change of variable $\eta=(ir^2\chi k b)^{1/3} (y-\lambda)$ gives the Airy equation as in section \ref{subsec:4.2}. Its solution can also be written in terms of the hypergeometric function $_0F_1$. We then obtain
\begin{align}
\Psi(y)&=d_1\mbox{Ai}[ (ir^2\chi kb)^{1/3} (y-\lambda)]+d_2\mbox{Bi}[ (ir^2\chi kb)^{1/3} (y-\lambda)] 
\nonumber\\&=d_1\ _0 F_1\left(;\frac{2}{3}; \frac{ir^2\chi kb(y-\lambda)^3}{9}\right)+d_2(y-\lambda)\ _0 F_1\left(;\frac{4}{3}; \frac{ir^2\chi kb(y-\lambda)^3}{9}\right), 
\label{eq86}
\end{align}
where $d_1$ and $d_2$ are constants. 
To approximate $\phi$, we integrate $\Psi$ twice and we obtain, see (\ref{eqC1}) and (\ref{eqC2}) in Appendix \ref{AppC},
\begin{multline}
\phi(y)=\int_y\Psi(y)dy^2=d_1(y-\lambda)^2\ _1 F_2\left(\frac{1}{3};\frac{4}{3},\frac{5}{3}; \frac{ir^2\chi kb(y-\lambda)^3}{9}\right)\\+d_2(y-\lambda)^3 \ _2 F_3\left(\frac{2}{3},1;\frac{4}{3},\frac{5}{3},2; \frac{ir^2\chi kb(y-\lambda)^3}{9}\right)+d_3y +d_4, r\to0^+,
\label{eq87}
\end{multline}
where $d_1$,  $d_2$,  $d_3$ and $d_4$ are constants.

We expect the second term to go to infinite as $y\to\infty$ since it resulted from integrating twice $\mbox{Bi}$, and it is known that $\mbox{Bi}\to\infty$ as $y\to\infty$ \cite{AB}. Therefore, $d_2=d_3=d_4=0$ in order $\phi$ to satisfy the boundary condition $\phi(0)=\phi(\infty)=0$. Hence, the eignfunctions are given by 
\begin{equation}
\phi_n(y)=[y-\lambda(k_n,\omega_n)]^2\ _1 F_2\left(\frac{1}{3};\frac{4}{3},\frac{5}{3}; \frac{ir^2\chi k_nb[y-\lambda(k_n,\omega_n)]^3}{9}\right),  n=0,1,2,\cdots.
\label{eq88}
\end{equation}
where $d_n$ are constants.  To compute the eigenvalues or to obtain the dispersion relation, we solve $\phi_y(y=0;k,\lambda(k,\omega))=0$, or equivalently, we find the zeros of  

\begin{equation}
 _1 F_2\left(\frac{1}{3};\frac{2}{3},\frac{4}{3}; \frac{-ir^2\chi kb\lambda^3}{9}\right)=0.
\label{eq89}
\end{equation}

\subsection{Configuration 4 ($r\to0^+$): $\bar{U}(y) = \delta y^2+by + c$, $\delta=$ a small constant}

We consider the mean flow profile given by $U(y) = \delta y^2 + by + c$, where $\delta$ is a small constant. We also consider that $U_{yy}=\delta\sim O(r^2)$ to make sure the third term in ({\ref{eq83}) is negligible compared to the other terms and can therefore be dropped. This gives

\begin{equation}
\phi_{yyyy}-[2r^2k^2+ir^2\chi k(\delta y^2+by+c-\omega/k)]\phi_{yy}=0.
\label{eq90}
\end{equation}
Setting $\Psi=\phi_{yy},\lambda=r\sqrt{i\chi k}[\omega/k+b^2/(2\delta)-c-2k/(i\chi)]/\sqrt{\delta}$, making the change of variable $\eta=(i r^2\chi k \delta)^{1/4} [y+b/(2\delta)]$ gives Hermite equation as in section \ref{subsec:4.2}. We then obtain, in terms of the Reynolds number, 
\begin{equation}
\Psi_{n}\sim\frac{\exp\left(-(iR k_{_{n}})^{\frac{1}{2}} \left(y+\frac{b}{2\delta}\right)^2\right)}{\sqrt{2^n 2}\,\pi^{\frac{1}{4}}} H_n\left((iR k_n)^{\frac{1}{4}}\left(y+\frac{b}{2\delta}\right)\right), n=0,1,2,\cdots, 
\label{eq91}
\end{equation} 
where as before $H_n$ are Hermite polynomials of order $n$.

We note that we have to chose Hermite polynomials with even order in order the solution to satisfy the boundary condition $ \phi_y(0)=\phi_y(\infty)=0$, and then integrate (\ref{eq91}) twice to approximate the eigen-solutions for (\ref{eq83})-(\ref{eq84}). The dispersion relation can be approximated as $\lambda(k_{m},\omega_{m})=2m+({1}/{4}), m=0,1,2,\cdots$. In that case, in the laminar boundary layer, the eigenvalues can be approximated by 
\begin{equation}
k_m\sim\frac{\sqrt{3}+i}{2}\left(\frac{R}{4}\right)^{1/3}\left(2m+\frac{1}{4}\right)^{2/3}, m=0,1,2,\cdots,
\label{eq92}
\end{equation} 
while we have
\begin{equation}
\Psi_{m}\sim\frac{\exp\left(-(iR k_{_{2m}})^{\frac{1}{2}} \left(y+\frac{b}{2\delta}\right)^2\right)}{\sqrt{2^m 2}\,\pi^{\frac{1}{4}}} H_{2m}\left((iR k_{2m})^{\frac{1}{4}}\left(y+\frac{b}{2\delta}\right)\right), m=0,1,2,\cdots.
\label{eq93}
\end{equation} 
\emph{Example~1.}
If $m=0$, for example, $H_0(y)=1$. Then using Proposition 1 in Nijimbere \cite{N} gives
\begin{multline}
\phi_{0y}=\int\Psi_0dy=\frac{1}{2\pi^{\frac{1}{4}}}\int e^{-(iR k_0)^{\frac{1}{2}}\left(y+\frac{b}{2\delta}\right)^2} H_0(y) dy\\=\frac{\left(y+\frac{b}{2\delta}\right)}{2\pi^{\frac{1}{4}}}\ _1F_1\left(\frac{1}{2},\frac{3}{2};-(iR k_0)^{\frac{1}{2}}\left(y+\frac{b}{2\delta}\right)^2\right), 
\label{eq94}
\end{multline}
where the constant of integration is set to zero in order to satisfy the boundary condition $\phi_y(0)=\phi_y(\infty)=0$.
Integrating (\ref{eq94}), see equation (\ref{eqC3}) in \ref{AppC}, gives 
\begin{equation}
\phi_{0}=\int\psi_{0y}dy=\frac{\left(y+\frac{b}{2\delta}\right)^2}{4\pi^{\frac{1}{4}}}\ _2F_2\left(\frac{1}{2},1;\frac{3}{2},2;-(iR k_0)^{\frac{1}{2}}\left(y+\frac{b}{2\delta}\right)^2\right),
\label{eq95}
\end{equation}
where the constant of integration is set to zero in order to satisfy the boundary condition $\phi(0)=\phi(\infty)=0.$
\\[1ex]

\section{Concluding remarks and discussion}\label{sec:6}

Squire's transformation were used to transform the three-dimensional model into a two-dimensional one. In two dimensions, the mean flow became  $\bar{U}(y) = \bar{u}\cos\theta + \bar{w}\sin\theta$ where $\bar{u}$, $\bar{w}$ are the streamwise and spanewise background fluid flow respectively, and $\theta$ is the phase velocity orientation angle in the horizontal plane ($xz$-plane). For instance, if  $\theta=\pi/6$, $\bar{u}=2(y-\tanh{y})$, $\bar{w}=(2/\sqrt{3})\tanh{y}$, then  $\bar{U}(y) =2(y-\tanh{y})\cos(\pi/6) + (2/\sqrt{3})\tanh{y}\sin(\pi/6)=y$. This applies to $\bar{U}(y)=y^2$ showing that a large class of three-dimensional background mean flow profiles can be represented in two dimensions by $\bar{U}(y) = ay^2 + by + c,$ $a, b$ and $c$ being constants, using Squire's transformation \cite{SH}. For the two-dimensional wake, the mean flow profile, $\bar{U}(y)=U_0\  \mbox{sech}^2(\varphi y), -\infty<y<+\infty, U_0>0, \varphi>0$, having a form of a Bell (Gaussian) function was used.

Making non-dimensional all variables and parameters, a spacial aspect ratio was introduced. This mainly allowed us to consider two configurations, the short-wave limit approximation which is obtained by letting the aspect ratio taking large values, and the long-wave limit approximation in which the aspect ratio takes small values in the Orr-Sommerfeld equation. This also allowed us to utilize analytical and asymptotic methods to obtain asymptotic solutions of the Orr-Sommerfeld equation and their corresponding eigenvalues. Most importantly, the procedure used in the present paper works regardless of the value of the Reynolds number.

In the short-wave limit approximation, the dispersion relation, the asymptotic eigenvalues and their corresponding asymptotic eigenfunctions were derived for configurations where the velocity mean flow profiles can, using Squire's transformation, be represented either as linear function or as a quadratic function. The eigenvalues were approximated using WKB method. Asymptotic eigenvalues and asymptotic eigenfunctions were also derived for the two-dimensional wake. The eigenfunctions were written in terms of Green's functions, and their corresponding outer approximate solutions were obtained as well. The results showed that the amplitude of the wave become larger as the Reynolds number becomes small. This is in agreement with the fact that small viscosity induces viscous instabilities. In the long-wave limit approximation, solutions were derived in terms of hypergeometric functions whose properties are known.

Due to the evolution of computer technology, Computational Fluid Dynamics (CFD) should help more in predicting transition from laminar flows to turbulent flows in three-dimensional shear flows. But starting simulations within a good range of eigenvalues allowing convergence of simulations to correct solutions remains a challenging problem. For this reason, we have approximated eigenvalues in the present paper which may, for instance, be used as a starting point in CFD simulations.


\begin{figure}
\vspace*{3mm}
\centerline{\includegraphics*[height=80mm,width=80mm]{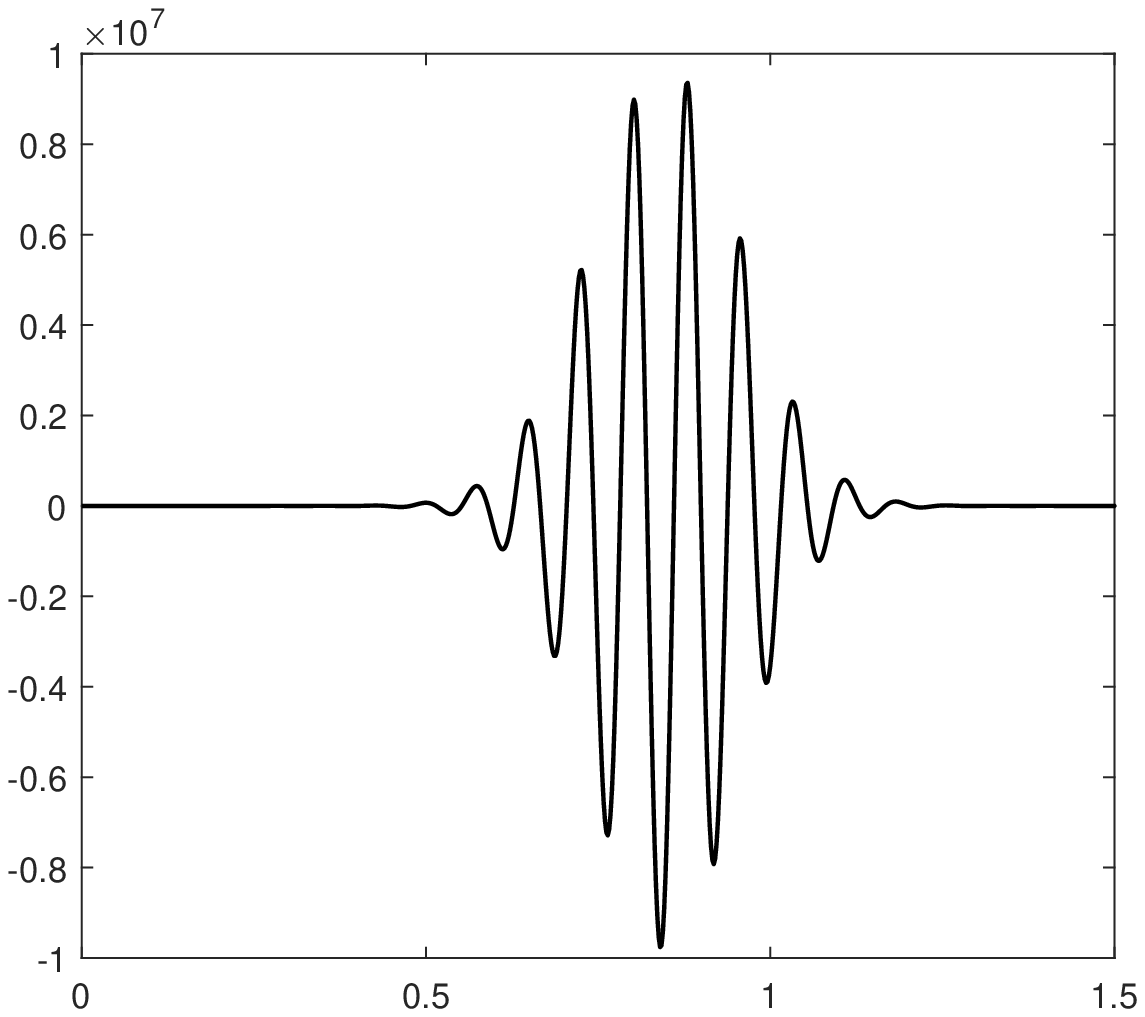}\includegraphics*[height=80mm,width=80mm]{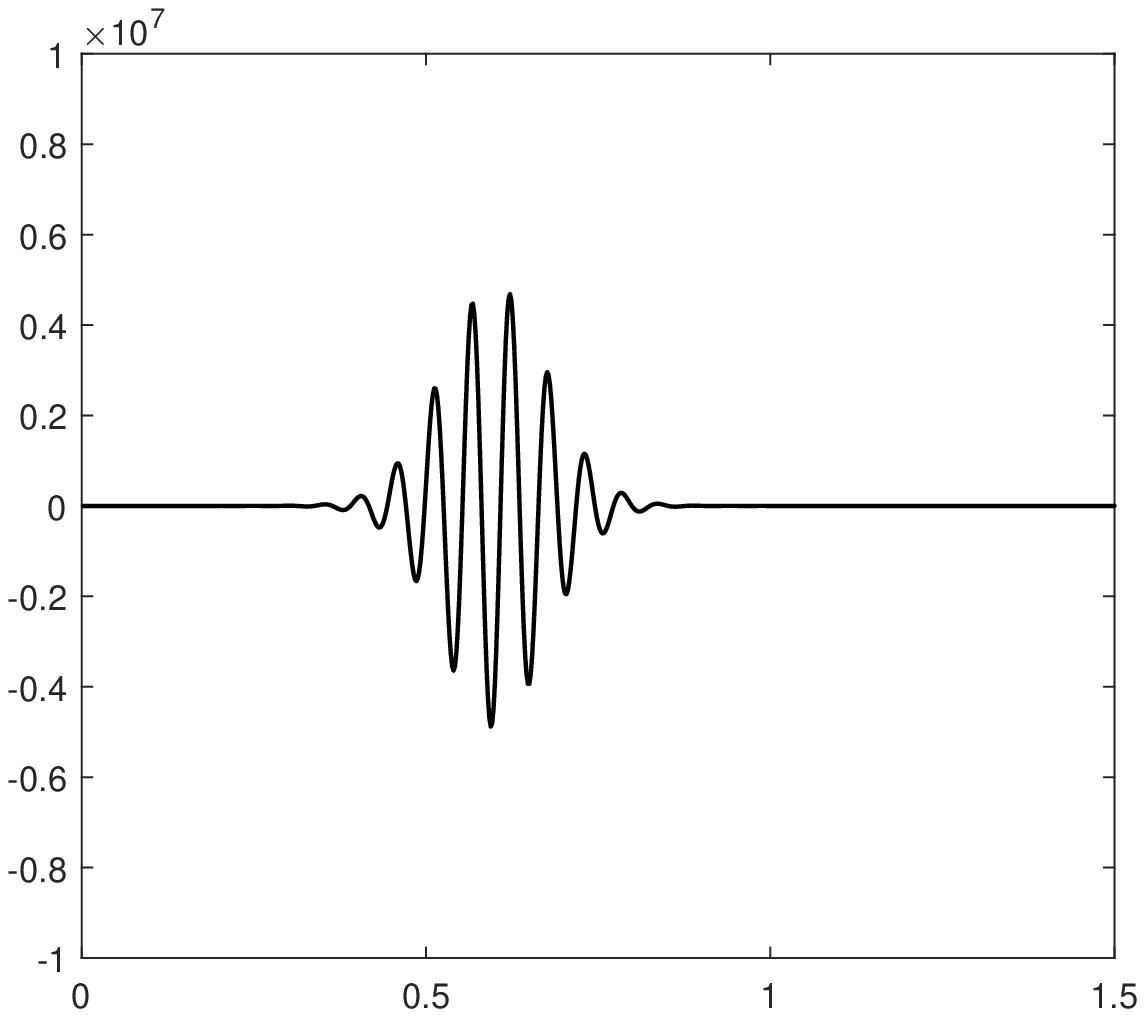}}
\vspace*{4mm}
\centerline{\includegraphics*[height=80mm,width=80mm]{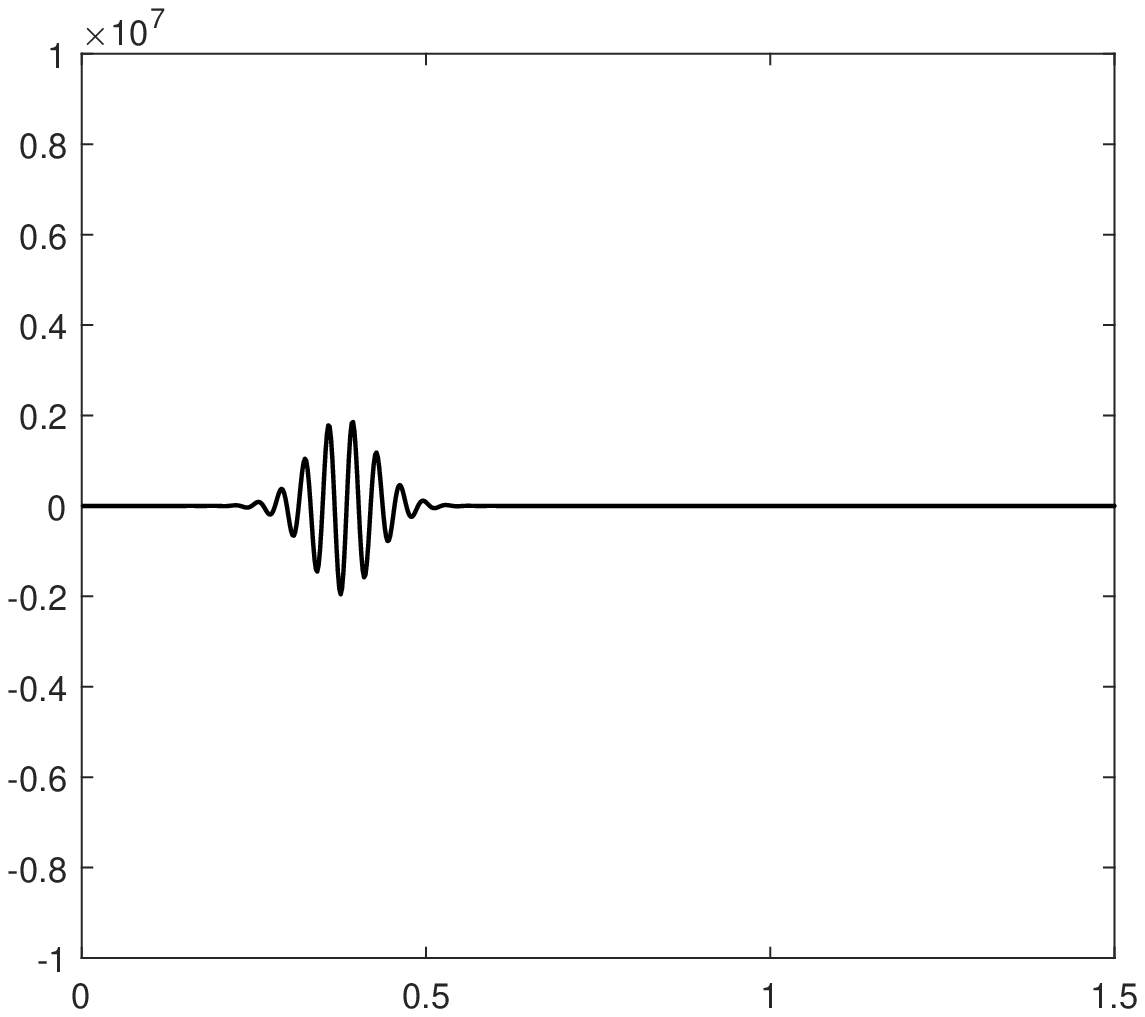}\includegraphics*[height=80mm,width=80mm]{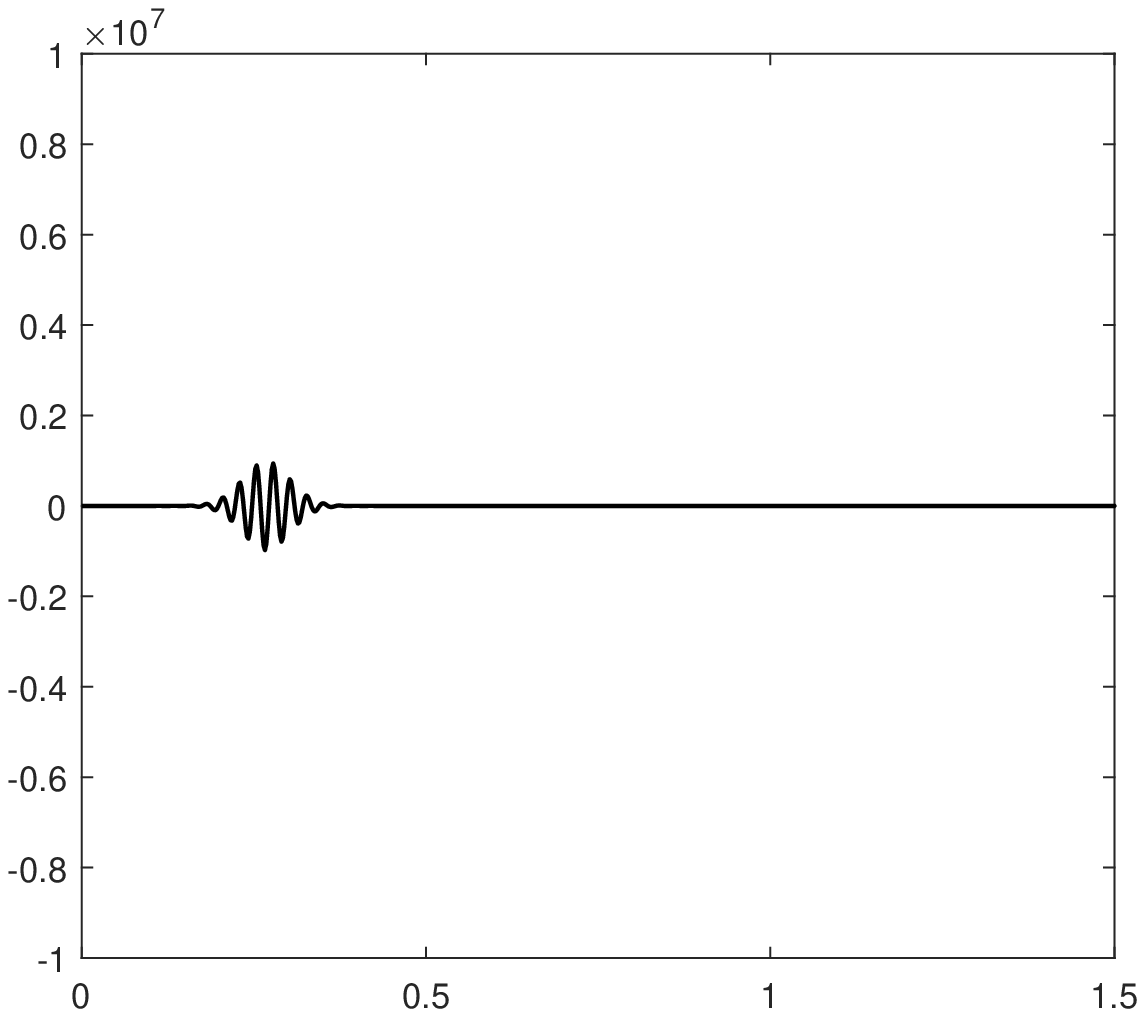}}
\begin{picture}(0,0)
\put (-35,370){\large $\phi_5$}
\put (190,370){\large $\phi_5$}
\put (-35,130){\large $\phi_5$}
\put (190,130){\large $\phi_5$}
\put (85,20){\large $y$}
\put (310,20){\large $y$}
\put (85,258){\large $y$}
\put (310,258){\large $y$}
\put (5,480){\large $(a)$ $R=1000$}
\put (225,480){\large $(b)$ $R=2000$}
\put (5,240){\large $(c)$ $R=5000$}
\put (225,240){\large $(d)$ $R=10000$}
\end{picture}
\caption{Plot of the eigenfunction $\phi_5$ as a function $y$. The corresponding Reynolds number are (a) $R=1000$, (b) $R=2000$, (c) $R=5000$ and (c) $R=10000$, $\bar{U}(y)=y, 0<y<\infty$ and the parameter $\epsilon=0.2$.}
\label{fig1}
\normalsize
\end{figure}

\begin{figure}
\centerline{\includegraphics*[height=80mm,width=80mm]{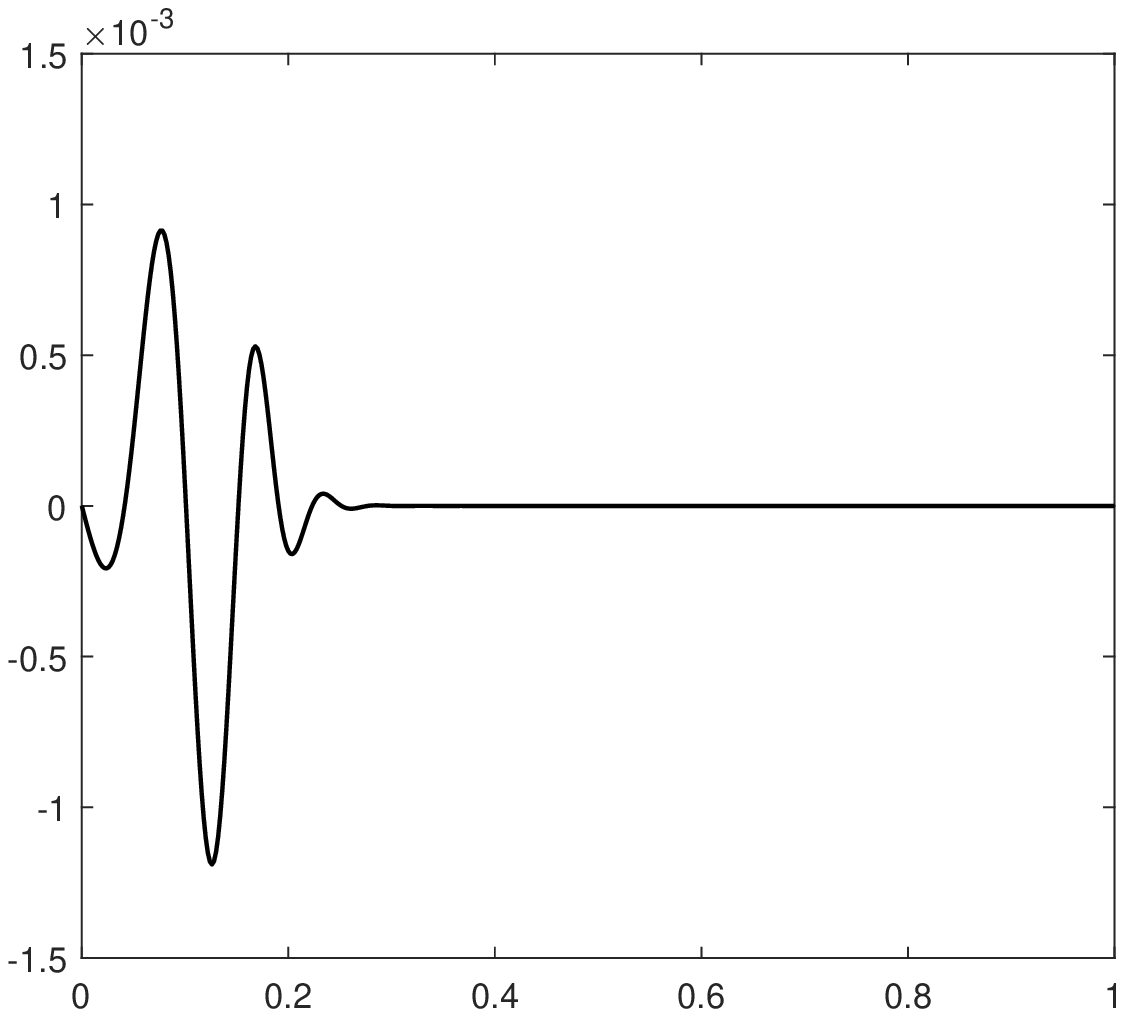}\includegraphics*[height=80mm,width=80mm]{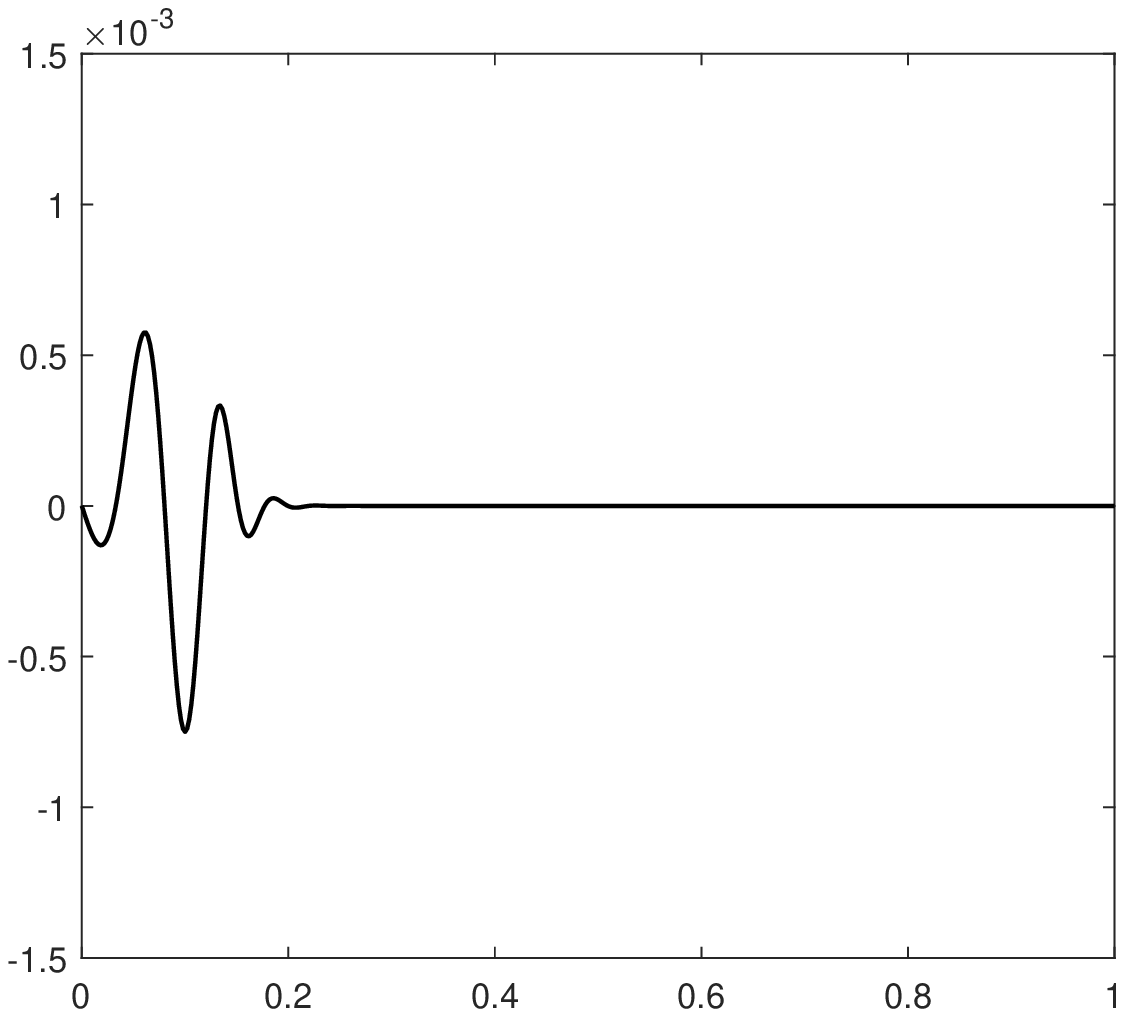}}
\vspace*{3mm}
\centerline{\includegraphics*[height=80mm,width=80mm]{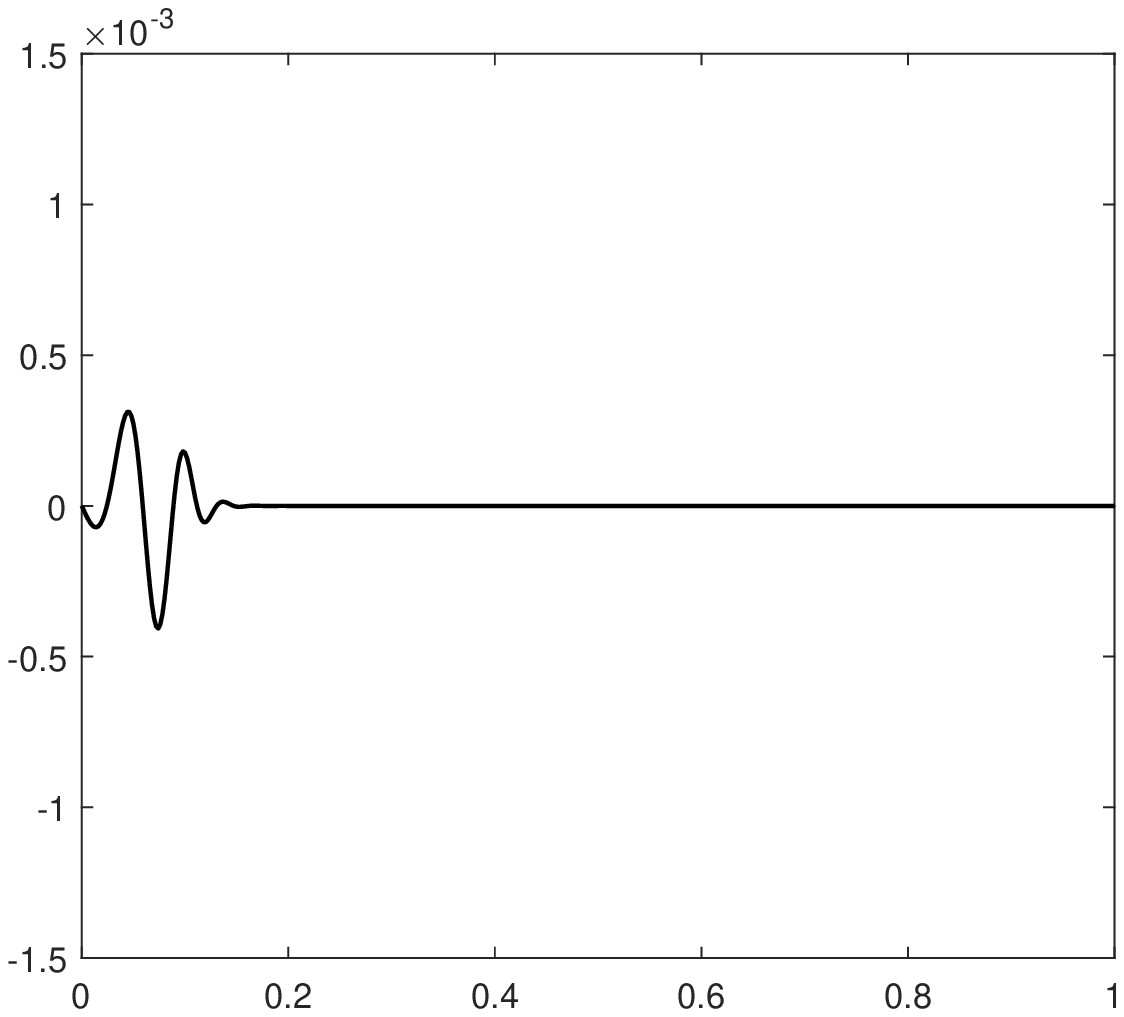}\includegraphics*[height=80mm,width=80mm]{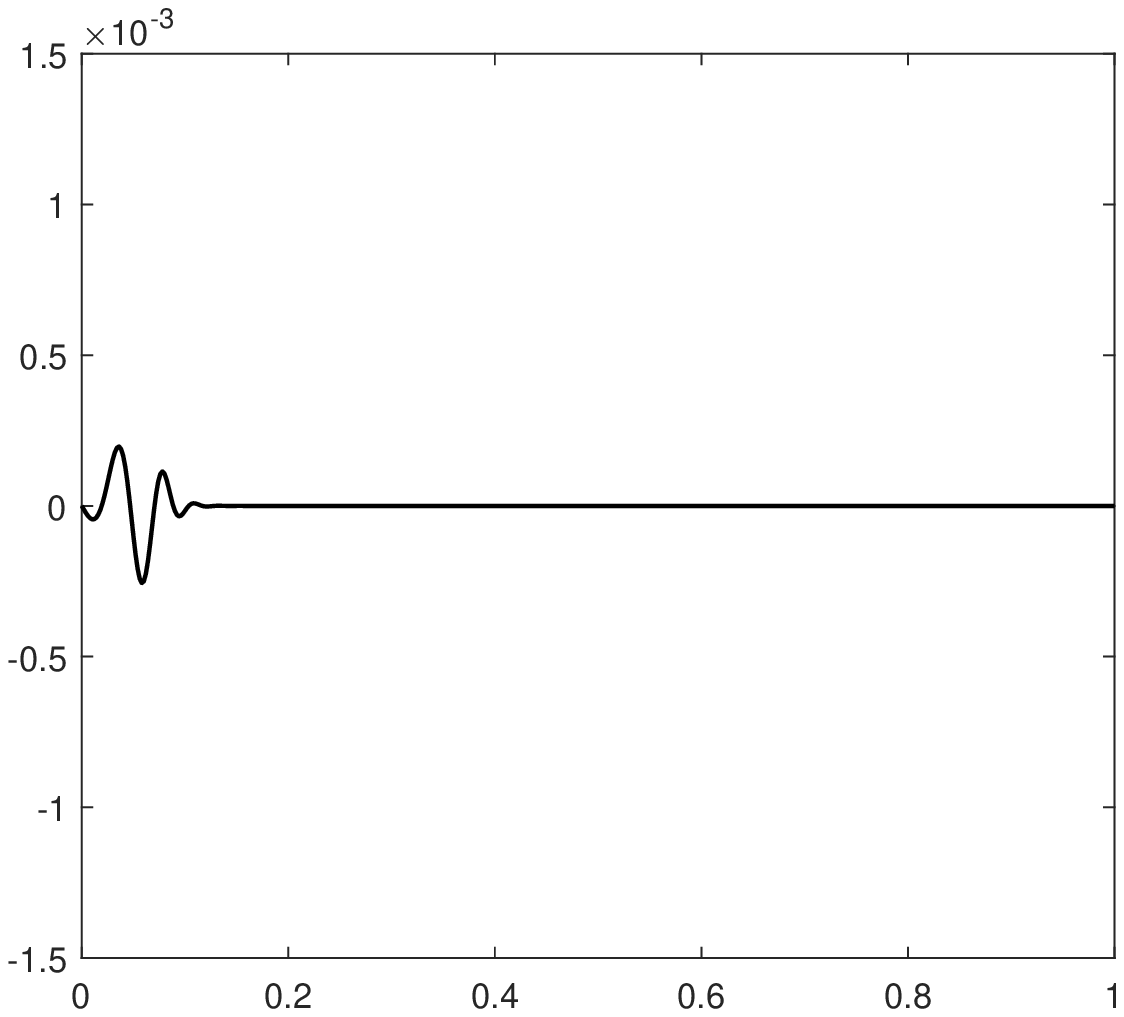}}
\begin{picture}(0,0)
\put (-35,370){\large $\phi_2$}
\put (190,370){\large $\phi_2$}
\put (-35,130){\large $\phi_2$}
\put (190,130){\large $\phi_2$}
\put (85,20){\large $y$}
\put (310,20){\large $y$}
\put (85,258){\large $y$}
\put (310,258){\large $y$}
\put (5,480){\large $(a)$ $R=1000$}
\put (225,480){\large $(b)$ $R=2000$}
\put (5,240){\large $(c)$ $R=5000$}
\put (225,240){\large $(d)$ $R=10000$}
\end{picture}
\caption{Plot of the eigen function $\phi_2$ as a function $y$. The corresponding Reynolds number are (a) $R=1000$, (b) $R=2000$, (c) $R=5000$ and (c) $R=10000$, $\bar{U}(y)=y^2, 0<y<\infty$ and the parameter $\epsilon=0.2$.}
\label{fig2}
\normalsize
\end{figure}

\appendix
\section{Green's functions}\label{AppA}
We want to solve, using Green's functions, the equation
\begin{equation}
\phi_{yy}-r^2k^2\phi=\Psi 
\label{eqA1}
\end{equation}
subject to boundary conditions
\begin{equation}
\phi_y(0)=\phi_y(\infty)=0.
\label{eqA2}
\end{equation}
In that case, Green's functions are some functions $G(y,\xi)$ that solve the equation 
\begin{equation}
G_{yy}-r^2k^2G=\delta(y-0) 
\label{eqA3}
\end{equation}
where $\delta$ is the delta Dirac function, and the functions $G$ satisfy the boundary conditions
\begin{equation}
G_y(0,\xi)=G_y(\infty,\xi)=0.
\label{eqA4}
\end{equation}

Green's functions are given by 
\begin{equation}
G(y,\xi)= \left\{\begin{array}{l}
\frac{u_1(y)u_2(\xi)}{W}, \mbox{if $y<\xi$ } \\
\frac{u_1(\xi)u_2(y)}{W}, \mbox{if $y>\xi$ }\\
 \end{array}\right., 
\label{eqA5}
\end{equation}
where $u_1$ and $u_2$ are solutions of the homogeneous equation 
\begin{equation}
u_{yy}-r^2k^2u=0 
\label{eqA6}
\end{equation}
subject to boundary conditions
\begin{equation}
u_y(0)=u_y(\infty)=0,
\label{eqA7}
\end{equation}
with $u_1$ satisfying the boundary conditions $u_{1y}(0)=0$, while $u_2$ satisfies $u_{2y}(\infty)=0$, and $W(y)=u_1u_{2y}-u_{1y}u_{2}$ is the associated Wronskian.
Then, $u_1(y)=e^{rky}+e^{-rky}$, $u_2(y)=e^{-rky}$, and so $W=-2rk$. This gives
\begin{equation}
G(y,\xi)= \left\{\begin{array}{l}
\frac{\cosh(rky)e^{-rk\xi}}{rk}, \mbox{if $y<\xi$ } \\
\frac{\cosh(rk\xi)e^{-rky}}{rk}, \mbox{if $y>\xi$ }\\
 \end{array}\right..
\label{eqA8}
\end{equation}

If on the other hand $\phi$ satisfies the boundary conditions 
\begin{equation}
\phi_y(-\infty)=\phi_y(\infty)=0, 
\label{eqA9}
\end{equation}
then Green's functions are given by 
\begin{equation}
G(y,\xi)= \left\{\begin{array}{l}
\frac{e^{rky}e^{-rk\xi}}{2rk}, \mbox{if $y<\xi$ } \\
\frac{e^{rk\xi}e^{-rky}}{2rk}, \mbox{if $y>\xi$ }\\
 \end{array}\right..
\label{eqA10}
\end{equation}

\section{WKB method approximation for the eigenvalues}\label{AppB}

We are interested in computing the eigenvalues for the Schrodinger equation
\begin{equation}
\epsilon u_{yy}-[V(y)-E]u=0, \epsilon\to0 
\label{eqB1}
\end{equation}
subject to boundary conditions
\begin{equation}
u(0)=u(\infty)=0. 
\label{eqB2}
\end{equation}
The physical-optics approximation requires that the valid solution in the region between $(0,\infty)$ is a linear combination of two rapidly oscillating WKB expressions \cite{BO}
\begin{equation}
u(y)=A [E-V(y)]^{-1/4}e^{\frac{i}{\epsilon}\int\limits_0^y\sqrt{E-V(\mu)} d\mu}+B [E-V(y)]^{-1/4}e^{-\frac{i}{\epsilon}\int\limits_0^y\sqrt{E-V(\mu)} d\mu}, 
\label{eqB3}
\end{equation}
where $A$  and $B$ are constants.

If $b\in(0,\infty)$, and 0 and $b$ are the turning points of $P$, i.e. $E-V(0)=0$ and $E-V(b)=0$, then the constants $A$ and $B$ must be chosen so that \cite{BO}
\begin{equation}
u(y)=2C [E-V(y)]^{-1/4}\sin\left[{\frac{1}{\epsilon}\int\limits_0^y\sqrt{E-V(\mu)} d\mu}+\frac{\pi}{4}\right]. 
\label{eqB4}
\end{equation}
Therefore, the eigenvalues have to satisfy 
\begin{equation}
\frac{1}{\epsilon}\int\limits_0^y\sqrt{E-V(\mu)}d\mu+\frac{\pi}{4}=n\pi, \epsilon\to0, n=0,1,2,\cdots. 
\label{eqB5}
\end{equation}
Hence, we have
\begin{equation}
\int\limits_0^y\sqrt{E-V(\mu)}d\mu\sim\epsilon\left(n-\frac{1}{4}\right)\pi, \epsilon\to0, n=0,1,2,\cdots, 
\label{eqB6}
\end{equation}
or
\begin{equation}
\int\limits_0^y\sqrt{E-V(\mu)}d\mu\sim\left(n-\frac{1}{4}\right)\frac{\pi}{r}, r\gg1, n=0,1,2,\cdots. 
\label{eqB7}
\end{equation}
in terms of the aspect ratio.

\section{Some useful integrals}\label{AppC}
Some useful integrals involving the generalized hypergeometric function (\ref{eq73}) are evaluated here. The method used is similar to that in Nijimbere\cite{N}.

\begin{align}
\int \ _0F_1\left(;\frac{2}{3};\frac{y^3}{9}\right) dy^2&=\int\sum\limits_{n=0}^\infty\frac{\left(\frac{y^3}{9}\right)^n}{\left(\frac{2}{3}\right)_n n!}dy^2=\sum\limits_{n=0}^\infty\frac{\left(\frac{1}{9}\right)^n}{\left(\frac{2}{3}\right)_n n!}\int y^{3n}dy^2\nonumber \\&=\sum\limits_{n=0}^\infty\frac{\left(\frac{1}{9}\right)^n}{\left(\frac{2}{3}\right)_n n!} \frac{y^{3n+2}}{(3n+1)(3n+2)}+C_1 y+C_2
\nonumber \\&=\frac{\Gamma\left(\frac{2}{3}\right)y^2}{9}\sum\limits_{n=0}^\infty\frac{\left(\frac{1}{9}\right)^n}{\Gamma\left(n+\frac{2}{3}\right) n!} \frac{y^{3n}}{\left(n+\frac{1}{3}\right)\left(n+\frac{2}{3}\right)}+C_1 y+C_2
\nonumber \\&=\frac{\Gamma\left(\frac{2}{3}\right)y^2}{9}\sum\limits_{n=0}^\infty\frac{\Gamma\left(n+\frac{1}{3}\right) }{\Gamma\left(n+\frac{4}{3}\right)\Gamma\left(n+\frac{5}{3}\right)}\frac{\left(\frac{y^3}{9}\right)^n}{n!}+C_1 y+C_2
\nonumber \\&=y^2\sum\limits_{n=0}^\infty\frac{\left(\frac{1}{3}\right)_n }{\left(\frac{4}{3}\right)_n\left(\frac{5}{3}\right)_n}\frac{\left(\frac{y^3}{9}\right)^n}{n!}+C_1 y+C_2
\nonumber \\&=y^2\ _1F_2\left(\frac{1}{3};\frac{4}{3},\frac{5}{3};\frac{y^3}{9}\right)+C_1 y+C_2.
\label{eqC1}
\end{align}
\begin{align}
\int y \ _0F_1\left(;\frac{4}{3};\frac{y^3}{9}\right) dy^2&=\int y\sum\limits_{n=0}^\infty\frac{\left(\frac{y^3}{9}\right)^n}{\left(\frac{4}{3}\right)_n n!}dy^2=\sum\limits_{n=0}^\infty\frac{\left(\frac{1}{9}\right)^n}{\left(\frac{4}{3}\right)_n n!}\int y^{3n+1}dy^2\nonumber \\&=\sum\limits_{n=0}^\infty\frac{\left(\frac{1}{9}\right)^n}{\left(\frac{4}{3}\right)_n n!} \frac{y^{3n+3}}{(3n+2)(3n+3)}+C_1 y+C_2
\nonumber \\&=\frac{y^3}{9}\sum\limits_{n=0}^\infty\frac{\left(\frac{1}{9}\right)^n}{\left(\frac{4}{3}\right)_n n!} \frac{y^{3n}}{\left(n+\frac{2}{3}\right)\left(n+1\right)}+C_1 y+C_2
\nonumber \\&=\frac{y^3}{9}\sum\limits_{n=0}^\infty\frac{\Gamma\left(n+\frac{2}{3}\right) \Gamma\left(n+1\right)}{\left(\frac{4}{3}\right)_n\Gamma\left(n+\frac{5}{3}\right) \Gamma\left(n+2\right)}\frac{\left(\frac{y^3}{9}\right)^n}{n!}+C_1 y+C_2
\nonumber \\&=\frac{y^3}{6}\sum\limits_{n=0}^\infty\frac{\left(\frac{2}{3}\right)_n (1)_n }{\left(\frac{4}{3}\right)_n\left(\frac{5}{3}\right)_n(2)_n}\frac{\left(\frac{y^3}{9}\right)^n}{n!}+C_1 y+C_2
\nonumber \\&=\frac{y^3}{6}\ _2F_3\left(\frac{2}{3},1;\frac{4}{3},\frac{5}{3},2;\frac{y^3}{9}\right)+C_1 y+C_2.
\label{eqC2}
\end{align}
\begin{align}
\int y \ _1F_1\left(\frac{1}{2}; \frac{3}{2};\gamma y^2\right) dy&=\int y\sum\limits_{n=0}^\infty\frac{\left(\frac{1}{2}\right)_n(\gamma y^2)^n}{\left(\frac{3}{2}\right)_n n!}dy\nonumber \\&=\sum\limits_{n=0}^\infty\frac{\left(\frac{1}{2}\right)_n \gamma^n}{\left(\frac{3}{2}\right)_n n!}\int y^{2n+1}dy\nonumber \\&=\sum\limits_{n=0}^\infty\frac{\left(\frac{1}{2}\right)_n\gamma^n}{\left(\frac{3}{2}\right)_n n!} \frac{y^{2n+2}}{2n+2}+C
\nonumber \\&=\frac{y^2}{2}\sum\limits_{n=0}^\infty\frac{\left(\frac{1}{2}\right)_n\gamma^n}{\left(\frac{3}{2}\right)_n n!} \frac{y^{2n}}{\left(n+1\right)}+C
\nonumber \\&=\frac{y^2}{2}\sum\limits_{n=0}^\infty\frac{\left(\frac{1}{2}\right)_n\Gamma\left(n+1\right)}{\left(\frac{3}{2}\right)_n\Gamma\left(n+2\right)}\frac{\left(\gamma y^2\right)^n}{n!}+C\nonumber \\&=\frac{y^2}{2}\sum\limits_{n=0}^\infty\frac{\left(\frac{1}{2}\right)_n (1)_n }{\left(\frac{3}{2}\right)_n(2)_n}\frac{\left(\gamma y^2\right)^n}{n!}+C
\nonumber \\&=\frac{y^2}{2}\ _2F_2\left(\frac{1}{2},1;\frac{3}{2},2;\gamma y^2\right)+C
\label{eqC3}
\end{align}

\end{document}